\documentclass[12pt, a4paper]{amsart}

\usepackage{framed}

\usepackage{amsmath, amssymb, amsfonts, amsthm}
\usepackage[margin=2.0cm]{geometry}
\usepackage[]{setspace}
\usepackage[usenames,dvipsnames]{color}
\usepackage[,unicode]{hyperref}
\hypersetup{
    urlcolor=blue,
    colorlinks=true,
    linkcolor= blue,
    citecolor=blue,
}
\usepackage{comment}

\usepackage[all,matrix,arrow]{xy}
\usepackage{enumerate}
\usepackage{xcolor}
\usepackage[utf8]{inputenc}
\usepackage[T1]{fontenc}
\usepackage[OT2,T1]{fontenc}
\DeclareSymbolFont{cyrletters}{OT2}{wncyr}{m}{n}
\DeclareMathSymbol{\Sha}{\mathalpha}{cyrletters}{"58}
\usepackage{enumitem}
\setlist[enumerate]{label={$(\alph*)$}, itemsep=0ex,leftmargin=0.6cm,topsep=1ex}
\setlist[enumerate,2]{label={\bf\alph*}., itemsep=1ex,leftmargin=0.5cm,topsep=1ex}
\setlist[enumerate,3]{label={\bf\roman*}., itemsep=1ex,leftmargin=0.5cm,topsep=1ex}

\usepackage{bm}
\usepackage{libertine}
\usepackage[italic]{mathastext} 

\theoremstyle{plain}
\newtheorem{theorem}{Theorem}[section]

\newtheorem{corollary}[theorem]{Corollary}

\theoremstyle{definition}
\newtheorem{example}[theorem]{Example}
\newtheorem{definition}[theorem]{Definition}
\newtheorem{question}[theorem]{Question}

\theoremstyle{remark}
\newtheorem{remark}[theorem]{Remark}

\renewcommand\int{\mathbf{int}}

\numberwithin{equation}{section}

\DeclareMathOperator{\res}{res}

\title[]{{A survey of local-global methods for Hilbert's Tenth Problem}}
\author{Sylvy Anscombe}
\address{Sylvy Anscombe, Universit\'{e} Paris Cit\'{e} and Sorbonne Universit\'{e}, CNRS, IMJ-PRG, F-75013 Paris, France}
\email{sylvy.anscombe@imj-prg.fr}
\author{Valentijn Karemaker}
\address{Valentijn Karemaker, Utrecht University, Mathematical Institute, PO Box 80010, 3508 TA Utrecht, The Netherlands}
\email{V.Z.Karemaker@uu.nl}
\author{Zeynep Kisakürek}
\address{Zeynep Kisak\"urek, Mathematisches Institut der Heinrich-Heine-Universität Düsseldorf,
Building: 25.22., D-40225 Düsseldorf, Germany}
\email{kisakuer@uni-duesseldorf.de}
\author{Vlerë Mehmeti}
\address{Vler\"e Mehmeti, Sorbonne Université and Universit\'{e} Paris Cit\'{e}, CNRS, IMJ-PRG, F-75005 Paris, France}
\email{vlere.mehmeti@imj-prg.fr}
\author{Margherita Pagano}
\address{Margherita Pagano, Mathematisch Instituut, Niels Bohrweg 1, 2333 CA Leiden, the Netherlands}
\email{m.pagano@math.leidenuniv.nl}
\author{Laura Paladino}
\address{Laura Paladino, Universit\`a della Calabria, Ponte Bucci, Cubo 30B, 87036, Rende (CS), Italy}
\email{laura.paladino@unical.it}

\begin{document}
\pagestyle{plain}
\maketitle

\begin{abstract}
Hilbert's Tenth Problem 
({\bf H10})
for a ring $R$
asks for an algorithm
to decide correctly,
for each $f\in\mathbb{Z}[X_{1},\dots,X_{n}]$,
whether the diophantine equation
$f(X_{1},...,X_{n})=0$
has a solution in $R$.
The celebrated `Davis-Putnam-Robinson-Matiyasevich theorem'
shows that 
{\bf H10} for $\mathbb{Z}$ is unsolvable,
i.e.~there is no such algorithm.
Since then, Hilbert's Tenth Problem has been studied in a wide range of rings and fields.
Most importantly, for {number fields and in particular for $\mathbb{Q}$},
{\bf H10}
is still an unsolved problem.
Recent work of Eisenträger, Poonen, Koenigsmann, Park, Dittmann, Daans, and others, has dramatically pushed forward what is known in this area, and has made essential use of local-global principles for quadratic forms, and for central simple algebras.
We give a concise survey and introduction to this particular rich area of interaction between logic and number theory, without assuming a detailed background of either subject. We also sketch two further directions of future research, one inspired by model theory and one by arithmetic geometry.
\end{abstract}

\section{Introduction}

Hilbert's Tenth Problem 
({\bf H10})
for a ring $R$
asks for an algorithm
to decide correctly,
for each $f\in\mathbb{Z}[X_{1},\dots,X_{n}]$,
whether the diophantine equation
\[
f(X_{1},...,X_{n})=0
\]
has a solution in $R$.
A stronger form allows polynomials with coefficients in certain recursive subrings of $R$.

Davis, Putnam, and Robinson,
largely working in the 50s and 60s \cite{DPR},
approached {\bf H10} by studying diophantine sets:
a set $A\subseteq\mathbb{Z}^{n}$ is {\em diophantine}
if there is a polynomial
$f\in\mathbb{Z}[X_{1},\ldots,X_{m},Y_{1},\ldots,Y_{n}]$
such that
\begin{align*}
	(a_{1},\ldots,a_{n})\in A\Longleftrightarrow \exists x_{1}\ldots x_{m}\;f(x_{1},\ldots,x_{m},a_{1},\ldots,a_{n})=0,
\end{align*}
for all $a_{1},\ldots,a_{n}\in\mathbb{Z}$.
They proved that a negative solution to
{\bf H10} would follow from the existence of a diophantine subset of $\mathbb{Z}$ of
`exponential growth'.
Moreover they gave a negative solution to the related `exponential' version of {\bf H10} for $\mathbb{Z}$ in which $f$ is allowed to be an exponential polynomial.
In 1970,
Matiyasevich \cite{Mat}
completed the proof of the celebrated `DPRM theorem' by finding a diophantine set of exponential growth using the Fibonacci sequence:

\begin{theorem}[Davis-Putnam-Robinson-Matiyasevich] \label{DPRM}
{\bf H10} for $\mathbb{Z}$ is unsolvable,
i.e.~there is no algorithm that decides correctly whether or not each multivariable polynomial has a solution in $\mathbb{Z}$.
\end{theorem} 

Since then, Hilbert's Tenth Problem has been studied in a wide range of rings, fields, and associated structures, in particular valued fields. 
Indeed, the study of {\bf H10} in fields and rings
brings together a variety of results and techniques from model theory, including the theories of local fields,
due to Tarski \cite{Tarski}, Ax--Kochen \cite{AxKochen},
and Ershov \cite{Ershov}, especially uniform aspects of these theories, as well as definability in global fields, field arithmetic, and valuation theory.

In the case of function fields of positive characteristic,
continuing a long line of work beginning with Pheidas \cite{Phe91},
Eisentrager and Shlapentokh showed in~\cite{ES} that {\bf H10} -- in its stronger form -- is unsolvable if the field does not contain the algebraic
closure of a finite field; this in particular applies to global function fields. For results over function fields of characteristic zero, see e.g.~\cite{KimRoushC, KimRoushpadic, Eis07, Eis08, RS}.
In another direction, Rumely showed in \cite{Rum86} that {\bf H10} is solvable for the ring of algebraic integers.
For more comprehensive overviews of this subject, the reader is encouraged to consult
\cite{KoeSurvey,ShlBook}.

Perhaps the most important active area of research in this direction is related to {\bf H10} for $\mathbb{Q}$,
which is still an open problem.
Poonen showed in \cite{Poonen03} that, for certain
large sets of primes $S$, {\bf H10} over $\mathbb{Z}[S^{-1}]$ has a negative answer.
Another major breakthrough came through the use of local-global principles
to find new families of diophantine subsets of~$\mathbb{Q}$.
Firstly, in \cite{robinson_1949} Robinson used the
Hasse-Minkowski theorem, a local-global principle on quadratic forms, to show
that every valuation ring $\mathcal{O}$ of a number field~$K$ is diophantine over $K$,
and that $\mathbb{Z}$ is definable in $\mathbb{Q}$ by a first-order formula in the language of rings
--see Section~\ref{section:logic} for a brief introduction to the terminology of first-order logic.
Subsequent steps were made in
the work of Eisenträger \cite{Eis05}
who used the Albert-Brauer-Hasse-Noether theorem, a local-global principle for central simple algebras,
and of Poonen \cite{Poonen}.

The ideas of Eisentr\"ager and Poonen kickstarted a series of discoveries.
More precisely, in \cite{Poonen}, Poonen proved a striking theorem:
$\mathbb{Z}$ is `universal-existentially definable' in $\mathbb{Q}$; that is, it is described by a first-order formula that consists of universal quantifiers, followed by exisential quantifiers, followed by a quantifier-free formula.
Building on this result, Koenigsmann showed in \cite{Koe16} that $\mathbb{Q}\setminus\mathbb{Z}$ is diophantine in $\mathbb{Q}$; Park in \cite{Park} in turn generalized this from $\mathbb{Q}$ to general number fields. In each case, the authors make use of the Hasse-Minkowski local-global principle for quadratic forms.
In a new direction,
in \cite{Dit18}, Dittmann generalized this machinery from quaternion algebras to general central simple algebras, and replaced the Hasse-Minkowski theorem with the
Albert--Brauer--Hasse--Noether theorem.
This allowed him to show that
`irreducibility is diophantine',
or more precisely:
for each global field~$K$ and each $n\in\mathbb{N}$, the set
$\{ (a_{0},\ldots,a_{n-1})\in K^{n}\mid\text{$X^{n}+\sum_{i<n}a_{i}X^{i}$ is irreducible in $K[X]$} \}$
is diophantine over $K$.
Further use has been made of the approach using central simple algebras,
for example by Daans \cite{Daans21,Daans23} and Anscombe--Dittmann--Fehm~\cite{ADF2}.

The above shows this is a growing research line of rich interaction between logic and model theory, and algebra and number theory.
Moreover, these results suggest natural extensions and generalizations in several directions, which we explore in the last section of this survey: on the one hand, the diophantine sets studied in these papers and the local-global behaviour they exhibit can fit into a framework of \emph{diophantine families}. On the other hand, one may replace the Hasse-Minkowski and Albert-Brauer-Hasse-Noether theorems with other local-global principles arising naturally in arithmetic geometry, and obtain different diophantine sets.

{In this article, we survey the stream of work 
following \cite{robinson_1949,Eis05,Poonen}
in which local-global principles are used to great effect.
After a concise introduction to the topic, we survey how the Hasse-Minkowski and Albert-Brauer-Hasse-Noether theorems have been used to prove striking new definability results in the above-mentioned works of
Eisentr\"ager, Poonen, Koenigsman, Park, Dittmann and Daans.
In particular, in Subsection~\ref{ssec:comp} we compare the techniques and results of these papers, in which the authors have been concerned with similar questions.
We hope that such an analysis will be useful for future uses of local-global principles in this study of definability in number theory, as we discuss in Subsections~\ref{ssec:mod} and~\ref{ssec:LGP}.}\\

Roughly speaking, we say a local-global principle holds when a property is true over a global field if and only if it is true over all, or all but finitely many, completions of that global field. 
One of the main reasons why local-global principles can be useful in the study of arithmetic properties of global objects is that, typically, the analogous properties over complete local fields are more accessible.
The above-mentioned Hasse-Minkowski and Albert-Brauer-Hasse-Noether theorems are concerned with the property of a quadratic form being isotropic, resp.~of a central simple algebra splitting (or equivalently, a Brauer group element being trivial). 

In arithmetic geometry one is concerned with rational points and arithmetic invariants of algebraic varieties and may study local-global principles for the existence of rational points (see e.g.~\cite{Katz, Har_Sko, Par_Sur}) and for the behaviour of arithmetic invariants such as endomorphisms and automorphisms (see e.g.~\cite{Sut}, \cite{Vogt}). 

There are many other contexts where local-global questions can be
considered. Indeed, local-global principles are studied also in related research areas like arithmetic dynamics (see e.g.~\cite{Hsi_sil,Bar}), and recently, in arithmetic geometry over a base field which is not a global field but for instance a function field over a number field \cite{Kato}, or a function field over a discretely valued field  \cite{CHHKPS,CHHKPS_2,CHHKPS_3}. 
In the literature there are also examples of different kinds of local-global questions, not formulated in the setting of fields at all. In those cases only the basic idea of studying properties in relation to behaviour at primes is borrowed from the classical local-global setup.
For example, in  \cite{Sus} the local behaviour of polynomial rings $R[X]$, where $R$ is a commutative ring, is considered over
its localizations at maximal ideals. In \cite{Malle} 
some properties
of the characters of a group are encoded by those of the normalizers of its $p$-Sylow
subgroups.  
However, as stated above, in the present paper we will consider the
classical local-global background of a global field and its completions.\\

This survey is intended for a broad audience, so we do not assume a detailed background in either logic or number theory. 
Accordingly, its structure is as follows: in Section~\ref{section:logic} we briefly introduce the necessary tools and terminology from logic,
and likewise in Section~\ref{section:algebra} we introduce the algebraic tools, in particular quaternion algebras and central simple algebras. We define the notion of a local-global principle and give the Hasse-Minkowski and Albert-Brauer-Hasse-Noether theorems as two important examples.
In Section~\ref{section:survey} we survey the main results and methods contained in Robinson's seminal work and the above-mentioned papers \cite{Eis05,Poonen,Koe16,Park,Dit18,Daans21}, in which the authors adopt similar methods of proof, in particular by making use of the Hasse-Minkowski and Albert-Brauer-Hasse-Noether local-global principles. 
Finally, in Section~\ref{sec:LGPdefinability}, we make comparisons between the surveyed articles, in order to establish a common model-theoretic framework.
We explore two further directions: on the one hand, we describe one possible model-theoretic viewpoint of diophantine sets satisfying some local-global principle; while on the other hand, we propose two other local-global principles from arithmetic geometry that could be used to study diophantine sets.

\section{Prerequisites from model theory}
\label{section:logic}

We begin by introducing a few basic concepts from first-order logic.

\subsection{First-order logic and definable sets}
We give a brief introduction to the basic tools of mathematical logic: formulas, structures, and definable sets.
Generally speaking, logic is the study of
valid reasoning.

In the same spirit,
the sub-discipline of

mathematical logic is the study of valid symbolic reasoning
within formal languages,
and it includes applications to theories drawn from mathematics.
Prominent examples of formal languages are the so-called first-order languages,
in which a significant proportion of mathematics can be formalized.
These languages are of particular importance in model theory, which studies the 
relationship between first-order theories and the mathematical structures that model them.
In this survey we will restrict most of our discussion to the first-order language of rings.
 
For more information, especially for a more general perspective on model theory,
readers are encouraged to consult
\cite{Hodges,Mar02}.

The first-order {\em language of rings},
denoted $\mathfrak{L}_{\mathrm{ring}}$,
is constructed in several stages, as follows.
We begin with the \emph{symbols} of the language,
which are of two types:
the logical and the non-logical.

The logical symbols are 
\begin{itemize}
     \item[---]
     {\em connectives}:
	$ \wedge $ for conjunction, $ \vee $ for disjunction, $\neg$ for negation, $\rightarrow$ for implication, and $ \leftrightarrow$ for biconditional,
     \item[---]
     {\em variables}:
	 $x_{1}, x_{2}, x_{3}, \ldots,y,z,\ldots$,
     \item[---]
     {\em quantifiers}:
	$\forall$ for the universal quantifier {\em "for all"} and $\exists$ for the existential quantifier {\em "there exists"},
     \item[---]
     {\em punctuation}:
	$=$ for equality, and the parentheses $($ \; $)$.
     \end{itemize}

The non-logical symbols
are 
the {\em function symbols}
$+$, $\cdot$, and $-$,
and the {\em constant symbols} $0$, and~$1$.
All these symbols are combined according to the following recursive construction to form two families: the terms and the first-order formulas.
The {\em terms} are the smallest collection of finite strings in the above symbols containing all the variables, the constant symbols $0$ and $1$, and such that if
$t_{1}$ and $t_{2}$ are terms
then
$(t_{1}+t_{2})$,
$(t_{1}\cdot t_{2})$,
and
$(t_{1}-t_{2})$
are terms.
The (first-order) {\em formulas} are the smallest collection of finite strings in the above symbols containing
$t_{1}=t_{2}$, for all terms $t_{1}$, $t_{2}$,
and also containing
$(\neg\varphi)$,
$(\varphi\wedge\psi)$,
$(\varphi\vee\psi)$,
$(\varphi\rightarrow\psi)$,
$(\varphi\leftrightarrow\psi)$,
$(\forall x)\,\varphi$,
and
$(\exists x)\,\varphi$,
for all formulas $\varphi,\psi$ and each variable $x$.
Informally, a first-order formula is a well-formed `finitary' mathematical statement, where quantifiers may range over elements (as opposed to sets or other `higher order' objects).

Examples of first-order formulas in $\mathfrak{L}_{\mathrm{ring}}$ are:
\begin{align}
\label{eq:a} & ((1+1)+1) =0,\\
\label{eq:b} & (\forall x)(x \cdot y)=(y \cdot x),\\
\label{eq:c} & (\forall x) (\exists y)(y \cdot y) = x,\\
\label{eq:d} & (\forall x)(((x \cdot x)+(y\cdot y))= 1 \rightarrow (\neg(x \cdot y) = 0))
\end{align}
For readability we may make various informal simplifications.
For example, Equation~\eqref{eq:b} may be written as
$(\forall x) x\cdot y=y\cdot x$,
and \eqref{eq:c} may be written as
$(\forall x)(\exists y) y\cdot y=x$
without ambiguity.
In each of the examples \eqref{eq:b}--\eqref{eq:d} above, the variable $x$ is in the {\em scope} of a quantifier: in other words,
$x$ is {\em bound} by that quantifier.
If a variable is not bound by any quantifier, it is {\em free}.
For example, $y$ is free in \eqref{eq:b} and \eqref{eq:d}.
We write $\varphi(x_{1},\ldots,x_{n})$ for $\varphi$ to signify that the free variables of $\varphi$ are among
$x_{1},\ldots,x_{n}$.

An {\em $\mathfrak{L}_{\mathrm{ring}}$-structure}
$(R,+^{R},\cdot^{R},-^{R},0^{R},1^{R})$
is a set
$R$
together with
with binary operations
$+^{R}$,
$\cdot^{R}$,
and~$-^{R}$
on $R$,
and constants
$0^{R},1^{R}\in R$.
Thus every ring is naturally an $\mathfrak{L}_{\mathrm{ring}}$-structure, as we would hope,
but there are plenty of $\mathfrak{L}_{\mathrm{ring}}$-structures that are not rings.
When working only with rings, further notational simplifications are possible,
for example \eqref{eq:a} may be written as $1+1+1=0$
and
\eqref{eq:d} as
$\forall x\;(x^2 + y^2 = 1 \rightarrow xy \neq 0)$.
Moreover we usually write an $\mathfrak{L}_{\mathrm{ring}}$-structure simply as e.g.~$(R,+,\cdot,-,0,1)$,
without distinguishing between the structure and the language, when there is no risk of ambiguity.
We may even write simply $R$ for $(R,+,\cdot,-,0,1)$ when the operations and constants are clear from the context.

We define a notion of `truth':
given a ring $(R,+,\cdot,-,0,1)$ and elements $a_{1},\ldots,a_{n}\in~R$,
each formula
$\varphi(x_{1},\ldots,x_{n})$ may be interpreted
by substituting each element $a_{i}$ into the place of its corresponding variable $x_{i}$.
Then the truth of
$\varphi(a_{1},\ldots,a_{n})$
in $R$ may be determined from (really, {\em defined by}) the usual mathematical meanings of `and', `or', `implies', etc., and the universal and existential quantifiers.
We write $R\models\varphi(a_{1},\ldots,a_{n})$
to mean that $\varphi(a_{1},\ldots,a_{n})$ is true in $R$.
For example, the formula $\varphi(y)$ appearing in \eqref{eq:b} has one free variable, namely $y$;
and for a ring $R$ and an element $b\in R$,
$\varphi(b)$ is true in $R$ (i.e.,~$R\models\varphi(b)$) if and only if $b$ commutes multiplicatively with each element of $R$.
For a full definition of interpretation, and of the relation $\models$, see for instance \cite[Definition 1.1.6]{Mar02}.

A {\em sentence} is a formula with no free variables.
The truth of a sentence $\varphi$ in a structure $R$ does not depend on any elements $a_{i}$.
Thus, in a given structure, each sentence is either true or false, and is not both.
For example, Equations~\eqref{eq:a} and \eqref{eq:c} are sentences.
By contrast, a first-order formula
$\varphi(x_1, \ldots, x_n)$ with free variables $x_1, \ldots, x_n$
defines the set of tuples $(a_{1},\ldots,a_{n})$ such that $\varphi$ is true when each $x_{i}$ is replaced by $a_{i}$. 
More broadly:
A set $A\subseteq R^{n}$ is
\textit{definable}
if there is an $\mathfrak{L}_{\mathrm{ring}}$-formula
$\varphi(x_{1},\ldots,x_{n},z_{1},\ldots,z_{l})$
and there are {\em parameters} $c_{1},\ldots,c_{l}\in R$
such that
for all 
$(a_{1},\ldots,a_{n})\in R^{n}$
\begin{align*}
	(a_{1},\ldots,a_{n})\in A
&\Longleftrightarrow
	R\models\varphi(a_{1},\ldots,a_{n},c_{1},\ldots,c_{l}),
\end{align*}
i.e.~if and only if
$\varphi(a_{1},\ldots,a_{n},c_{1},\ldots,c_{l})$ is true in $R$.

There are several special kinds of definability,
distinguished by the syntactic complexity
of the defining formula, as follows.
A formula $\varphi$ is
\begin{itemize}
\item
{\em quantifier-free} if it contains no quantifiers;
\item
{\em diophantine} if it consists of existential quantifiers followed by a single polynomial equality with integer coefficients;
$f(x_{1},\ldots,x_{n},y_{1},\ldots,y_{m},z_{1},\ldots,z_{l})=0$;
\item
{\em existential} if it consists of existential quantifiers followed by a quantifier-free formula;
\item
{\em universal} if it consists of universal quantifiers followed by a quantifier-free formula;
\item
{\em universal-existential} if it consists of universal quantifiers followed by an existential formula.
\end{itemize}
For example, the formula~\eqref{eq:a} is quantifier-free, whereas the formulas~\eqref{eq:b} and \eqref{eq:d} are universal, and the formula~\eqref{eq:c} is universal-existential.

In the study of {\bf H10}, we are mostly concerned with diophantine and existential formulas and the sets they define.
Straight from the definition, every diophantine formula is existential.
Given a ring $R$ and 
a subset $A \subseteq R^n$,
we say $A$ is
{\em existentially} definable
(respectively, {\em diophantine})
if it is defined by an existential (respectively, diophantine) formula.
Likewise, $A$ is
{\em universally} (respectively, {\em universal-existentially}) definable
if it is defined by a 
universal (resp., universal-existential)
formula.
Existentially definable sets are preserved in extensions: if $R\subseteq S$ is an extension of rings, and $A\subseteq R^{n}$ and $B\subseteq S^{n}$ are both defined by the same existential formula,
then $A\subseteq B$ \cite[Theorem 2.4.1]{Hodges}.
Moreover,
one repeatedly uses that
finite intersections and finite unions of diophantine sets are diophantine.

\begin{remark}\label{rem:constants}
We might also consider a definition of `definability' for sets $A$ of a ring $R$ that requires that the parameters $c_{1},\ldots,c_{l}$ are elements of a given subring $C\subseteq R$.
We then say that $A$ is {\em definable over $C$}, and likewise for `existentially definable', etc.
In the majority of this article we make no such restriction, and for clarity we often write e.g.~`$A$ is diophantine over $R$' to remind the reader that both $A$ is definable in the structure $R$ and that the parameters may be any elements of $R$.
\end{remark}

\begin{remark}\label{rem:conjunction_elimination}
In many of the rings and fields we consider,
a set is diophantine if and only if it is existentially definable.
More precisely, there are several standard methods that may be used to 
find `reductions' of a given existential formula to an equivalent diophantine one.
These reductions work in various settings:
for example
if $R$ is a field then
\[R\models\forall x_{1}\ldots\forall x_{n}\;(f(x_{1},\ldots,x_{n})\neq 0\leftrightarrow \exists y\;f(x_{1},\ldots,x_{n})\cdot y=1),\]
so in fields we may eliminate a negation in favour of an extra existential quantifier.
If $R$ is an integral domain then 
\[\begin{array}{ll}
R\models &\forall x_{1}\ldots\forall x_{n}\;((f(x_{1},\ldots,x_{n})=0\vee g(x_{1},\ldots,x_{n})=0)\leftrightarrow \\
& f(x_{1},\ldots,x_{n})\cdot g(x_{1},\ldots,x_{n})=0),\\
\end{array}\]
for $f,g\in R[x_{1},\ldots,x_{n}]$.
Thus in integral domains disjunctions of polynomial equalities may be replaced by a single polynomial equality.
On the other hand, if $R$ is a ring that admits a total order (compatible with the addition and multiplication),
then $a^{2}+b^{2}=0$ if and only if $a=b=0$, for all $a,b\in R$.
Thus we have
\[\begin{array}{ll}
R\models & \forall x_{1}\ldots\forall x_{n}\; (f(x_{1},\ldots,x_{n})=g(x_{1},\ldots,x_{n})=0
\leftrightarrow \\ &(f(x_{1},\ldots,x_{n}))^{2}+(g(x_{1},\ldots,x_{n}))^{2}=0).\\
\end{array}\]
Therefore in such a ring a conjunction of polynomial equalities may be replaced by a single polynomial equality.
A similar method works whenever $R$ is an integral domain and there is a non-constant polynomial over $R$ with no zero in the field of fractions of $R$.
\end{remark}

\begin{example}{Examples} \hspace{0.1cm}  
\begin{enumerate}
\item
$\mathbb{R}\setminus\{0\}$ is diophantine over $\mathbb{R}$,
defined by the formula
$(\exists y)\;xy=1$.
\item
$\mathbb{R}_{\geq0}$ is diophantine over $\mathbb{R}$,
defined by the formula
$(\exists y)\;x=y^{2}$.
\item
$\mathbb{Q}_{\geq0}$ is diophantine over $\mathbb{Q}$,
defined by the formula
$(\exists w_{1})(\exists w_{2})(\exists w_{3})(\exists w_{4})\;x=w_{1}^{2}+w_{2}^{2}+w_{3}^{2}+w_{4}^{2}$,
which makes use of Euler's theorem on sums of rational squares.
\item
$\mathbb{N}$ is diophantine over $\mathbb{Z}$,
defined by the formula
$(\exists w_{1})(\exists w_{2})(\exists w_{3})(\exists w_{4})\;x=w_{1}^{2}+w_{2}^{2}+w_{3}^{2}+w_{4}^{2}$.
\item
$\mathbb{Z}\setminus\{0\}$ is diophantine over $\mathbb{Z}$,
defined by the formula
\[(\exists w_{1})(\exists w_{2})(\exists w_{3})(\exists w_{4})\;(w_{1}^{2}+w_{2}^{2}+w_{3}^{2}+w_{4}^{2}+1-x)(w_{1}^{2}+w_{2}^{2}+w_{3}^{2}+w_{4}^{2}+1+x) =0.
\]
\item
$(-1,1)$ is diophantine over $\mathbb{R}$,
defined by the formula
\[(\exists w)\; (1-x^2)w^2=1.\]
\end{enumerate}
\end{example}

\begin{example}

In \cite{robinson_1965}, Robinson
 showed that $\mathbb{Z}_{p}$ is diophantine over $\mathbb{Q}_{p}$:
for any prime $l\neq p$,
the diophantine formula
$\exists y\;1+px^{l}=y^{l}$
defines $\mathbb{Z}_{p}$ in $\mathbb{Q}_{p}$.
To see this:
if
$a,b\in\mathbb{Q}_{p}$ are such that
$1+pa^{l}=b^{l}$,
then a simple calculation of $p$-adic valuations shows that
$a$ lies in $\mathbb{Z}_{p}$.
Conversely, if $a\in\mathbb{Z}_{p}$
then $1+pa^{l}$ lies in $1+p\mathbb{Z}_{p}$,
and by Hensel's Lemma every such element is an $l$-th power.

\end{example}

In Section~\ref{Robinson} we will recall Robinson's definition of $\mathbb{Z}$ in $\mathbb{Q}$.

\subsection{Theories}\label{theories}
 
Remaining in the language of rings,
a {\em theory} is simply a set of sentences.
Accordingly, for a ring~$R$, we define the {\em theory of $R$},
denoted $\mathrm{Th}(R)$,
to be the set of all $\mathfrak{L}_{\mathrm{ring}}$-sentences that hold in $R$.
Similarly the {\em existential theory} of $R$,
denoted $\mathrm{Th}_{\exists}(R)$,
is the set of existential $\mathfrak{L}_{\mathrm{ring}}$-sentences that are true in $R$.
A theory $T$ {\em axiomatizes}
the class of its models, i.e.,~those $\mathfrak{L}_{\mathrm{ring}}$-structures $R$ such that $R\models\varphi$ for every $\varphi\in T$.
For example the class of fields is axiomatized by the usual axioms for fields: associativity and commutativity of additions, etc., including the axiom 
\[ (\forall x) (\exists y)\;(x=0\vee xy=1) 
\]
that describes the existence of multiplicative inverses.

Given a theory $T$, naturally we ask ourselves if there is an algorithm (i.e., a Turing machine)
that correctly determines whether or not any given sentence in the language of rings is true in all models of $T$.
When there is such an algorithm, we say that $T$ is \textit{decidable}; otherwise it is \textit{undecidable}.
For example, 
as a consequence of Gödel's famous First Incompleteness theorem,
the theory $\mathrm{Th}(\mathbb{Z})$ of the ring $\mathbb{Z}$ is undecidable.
In particular, Peano's axioms are incomplete.
This is in contrast with the situation for $\mathbb{R}$ and $\mathbb{C}$ whose complete theories are decidable and admit recursive axiomatizations, due to Tarski \cite{Tarski},
and for
$\mathbb{Q}_{p}$,
due to Ax--Kochen \cite{AxKochen}, and independently by Ershov \cite{Ershov}.

Robinson's first-order definition of $\mathbb{Z}$ in $\mathbb{Q}$
(\cite{robinson_1949}, and see Theorem~\ref{Rob_Thm_3.1})
proves us with an {\em interpretation} 
of the theory of $\mathbb{Z}$ inside the theory of $\mathbb{Q}$.
More precisely,
let $ \psi(x)$ be a formula such that for any $a \in \mathbb{Q}$,
\[ \mathbb{Q}\models \psi(a) \Longleftrightarrow a \in \mathbb{Z}.\]
Then we may translate each formula $\varphi(x)$
into a formula $\hat{\varphi}(x)$ so that
\begin{align}\label{eq:interpretation}
	\mathbb{Z}\models\varphi(a_{1},\ldots,a_{n})\Longleftrightarrow\mathbb{Q}\models\hat{\varphi}(a_{1},\ldots,a_{n}),
\end{align}
for any formula $\varphi(x)$ and any integers $a_{1},\ldots,a_{n}$.
We proceed by recursion on the recursive definition of formulas:
the terms and atomic formulas remain unchanged,
the connectives also are unaltered,
but we `relativize' quantifiers,
as follows.
Since
\[\mathbb{Z}\models \exists x  \; \varphi(x,b_{1},\ldots,b_{n}) \Longleftrightarrow \mathbb{Q}\models\exists x \;(\psi(x) \wedge \hat{\varphi}(x,b_{1},\ldots,b_{n}))\]
and
\[ \mathbb{Z}\models \forall x  \; \varphi(x,b_{1},\ldots,b_{n}) \Longleftrightarrow \mathbb{Q}\models\forall x \;(\psi(x) \rightarrow \hat{\varphi}(x,b_{1},\ldots,b_{n})),\]
for all tuples $b_{1},\ldots,b_{n}$ of integers,
we define
$\widehat{\exists x\;\varphi}(y_{1},\ldots,y_{n})$ to be
$\exists x\;(\psi(x)\wedge\hat{\varphi}(x,y_{1},\ldots,y_{n}))$
and
$\widehat{\forall x\;\varphi}(y_{1},\ldots,y_{n})$ to be
$\forall x\;(\psi(x)\rightarrow\hat{\varphi}(x,y_{1},\ldots,y_{n}))$.
This construction yields a primitive recursive function
$\varphi(y_{1},\ldots,y_{n})\mapsto\hat{\varphi}(y_{1},\ldots,y_{n})$,
from the set of $\mathfrak{L}_{\mathrm{ring}}$-formulas to itself
satisfying
\eqref{eq:interpretation}.
Thus for every $\mathfrak{L}_{\mathrm{ring}}$-sentence $\varphi$,
to check whether or not
$\mathbb{Z}\models\varphi$
it suffices to check whether or not
$\mathbb{Q}\models\hat{\varphi}$.
By Gödel's First Incompleteness theorem the theory of $\mathbb{Z}$ is undecidable,
and therefore the theory of $\mathbb{Q}$ is undecidable

Let us return to Hilbert's Tenth Problem ({\bf H10}) for $\mathbb{Q}$.
When we translate it into the language of mathematical logic, this is the question of whether the existential theory
$\mathrm{Th}_{\exists}(\mathbb{Q})$
is decidable. 
If $\mathbb{Z}$ is a diophantine subset of $\mathbb{Q}$,
then the interpretation 
$\varphi(x)\mapsto\widehat{\varphi}(x)$
restricts to a map from
existential sentences to existential sentences,
still satisfying~\eqref{eq:interpretation}.
Consequently,
by Theorem~\ref{DPRM},
if~$\mathbb{Z}$ is a diophantine subset of $\mathbb{Q}$ 
then ({\bf H10}) is also unsolvable for $\mathbb{Q}$.
This motivates the search for a diophantine definition of $\mathbb{Z}$ in $\mathbb{Q}$.

Koenigsmann's definition of $\mathbb{Z}$ in $\mathbb{Q}$
\cite[Theorem 1]{Koe16}
is by a universal formula.
Using 
this formula, 
the interpretation
$\varphi(x)\mapsto\hat{\varphi}(x)$
restricts to a map from existential sentences to universal-existential sentences,
and so that at least the `$\forall\exists$-theory' of $\mathbb{Q}$ is undecidable,
i.e., there is no algorithm that decides correctly whether or not a universal-existential $\mathfrak{L}_{\mathrm{ring}}$-sentence is true in $\mathbb{Q}$ \cite[Corollary~3]{Koe16}.

\section{Prerequisites from algebra}
\label{section:algebra}
In this section we first give a brief description of quaternion algebras and, more generally,
central simple algebras. 
More details can be found for example in \cite[Chapter 2]{GS06}.
Afterwards we recall some local-global principles relevant to work on Hilbert's Tenth Problem. 

\subsection{Quaternions and Central Simple Algebras}

Throughout, we let $K$ be a field; we will assume for simplicity that $\mathrm{char}(K) \neq 2$.
All $K$-algebras are assumed to be associative and to have a multiplicative identity.
Moreover we consider only $K$-algebras of finite dimension.

\begin{definition}
A $K$-algebra $A$ is \emph{simple} if it has no non-trivial two-sided ideals.
It is a \emph{central simple algebra} if moreover its centre is $K$.
\end{definition}

The dimension of a central simple algebra $A/K$ is a square integer, and its positive square root is called its {\em degree}. 
A central simple algebra $A/K$ is said to be \emph{split} if it is isomorphic to a matrix algebra $M_n(K)$ for some~$n$. Otherwise we call it \emph{non-split}.
For any field extension $L/K$, we may consider $A \otimes_K L$, which is a central simple $L$-algebra. 
For $K$ a global field and $\mathfrak{p}$ a place of $K$,
we will want to consider the $K_{\mathfrak{p}}$-algebra $A \otimes_K K_{\mathfrak{p}}$,
where $K_{\mathfrak{p}}$ denotes the completion of $K$ with respect to $\mathfrak{p}$.
The set $\Delta_{A/K}$ of places $\mathfrak{p}$ of $K$ at which $A \otimes_K K_{\mathfrak{p}}$ is non-split is always finite (we will also denote it simply by
$\Delta_{A}$, when $K$ is clear from the context).
 
\begin{definition} \label{quaternion}
For $a,b \in K^{\times}$, we denote by $H_{a,b}$ the \emph{quaternion algebra} over $K$ with generators $i,j$ satisfying $i^{2}=a$, $j^{2}=b$, and $ij=-ji$.
Thus, 
\begin{align*}
    H_{a,b}&=\big\{w+xi+yj+zij\;\big|\;w,x,y,z\in K\big\}.
\end{align*}
\end{definition}

Every four-dimensional central simple algebra over $K$ is isomorphic to $H_{a,b}$ for some $a, b \in K^\times$ \cite[Proposition 1]{EGA}.
Every quaternion algebra $H_{a,b}$ is either split, so isomorphic to the matrix algebra $M_{2}(K)$;
or it is a division algebra, in which case it is non-split.
For an extension $L/K$, we denote $H_{a,b}(L) := H_{a,b}\otimes_{K}L$. For a global field $K$ and a place $\mathfrak{p}$ of $K$, we say $H_{a,b}$ is \emph{split} at $\mathfrak{p}$ if $H_{a,b}(K_{\mathfrak{p}}) \simeq M_2(K_{\mathfrak{p}})$,
and is \emph{ramified} at $\mathfrak{p}$ otherwise.

\begin{definition} \label{trnr}
\begin{enumerate}
\item The \emph{reduced trace} of an element $\alpha = w+xi+yj+zij$ of $H_{a,b}$ is defined as $\mathrm{tr}(\alpha) = 2w$. 
\item Its \emph{reduced norm} is defined as $\mathrm{nr}(\alpha) = w^2 - ax^2 - by^2 + abz^2$.
\end{enumerate}

Any central simple algebra $A/K$ has reduced trace and norm maps,
defined on $A$ with values in~$K$ (cf.~\cite[\S 2.6]{GS06}), which we will also denote by $\mathrm{tr}$ and $\mathrm{nr}$, respectively. The reduced trace is $K$-linear and the reduced norm is multiplicative on $A^\times$.

Let $A$ be a degree $n$ central simple $K$-algebra.
For $\alpha \in A$, let $P_{\alpha}(X)=X^n+ a_{n-1}X^{n-1}+\cdots + a_0$ be the reduced characteristic polynomial of $\alpha$ over $K$. Then $\mathrm{tr}(\alpha)=a_{n-1}$ and $\mathrm{nr}(\alpha)=(-1)^na_0.$
\end{definition}

\begin{definition}\label{def: Brauer group}  
We say that two central simple algebras $A$ and $A'$ over $K$ are \emph{Brauer equivalent} if there exist integers $m,m'>0$ such that 
    \[
        A\otimes_K M_{m}(K)\simeq A'\otimes_K M_{m'}(K).
    \]
    This defines an equivalence relation on the set of central simple algebras over $K$, equipped with the tensor product. The set of equivalence classes is an abelian group denoted by $\mathrm{Br}(K)$, see \cite{CoSko} and \cite{Sko} for more details. Finally, the class of a central simple algebra $A/K$ is trivial in $\mathrm{Br}(K)$ if and only $A$ is split. 
\end{definition}

\subsection{Local-global principles}
Let $K$ be a global field and $\mathfrak{p}$ a place of $K$.
Let $M_K$ be the set of all places of $K$. 
A quadratic form $F(X_1,... ,X_n)\in
 K[X_1,... ,X_n]$ is \emph{isotropic} over $K$ if the equation $F=0$ has a nontrivial solution
 in $K$.

\begin{theorem}{\cite[Theorem VI.3.1]{Lam05}} \label{HM}
Let $F(X_1,... ,X_n)\in
 K[X_1,... ,X_n]$ be a quadratic form, where $K$ is a global field. If $F$ is isotropic over $K_{\mathfrak{p}}$, for all completions
$K_{\mathfrak{p}}$ of $K$, then $F$ is isotropic over $K$.
\end{theorem}

This result was proved by Minkowski in 1923 in the case when $K=\mathbb{Q}$
and by Hasse in 1924 in the case when $K$ is a number field.
Since then many mathematicians have been inspired to study similar so-called \emph{local-global problems}, i.e.,\ problems
in which the validity of a certain property is assumed to hold in all, or all but finitely many, completions of a global field
$K$ and one asks if the same property holds in $K$.  When the answer to a local-global problem is affirmative one says that
a local-global principle or a Hasse principle holds. 

The validity of certain local-global principles implies
properties of the global fields over which they are formulated. For instance, if in Theorem
\ref{HM} we take $F=X^2-a$ and $K=\mathbb{Q}$, then we get that a rational number $a$ is
a square modulo all primes $\mathfrak{p}$ if and only if it is a square in  $\mathbb{Q}$.
Thus we have $\mathbb{Q}^2=\bigcap_{\mathfrak{p}} \mathbb{Q}_\mathfrak{p}^2$.  
Hence, knowing the validity of a 
local-global principle allows us to write certain subsets of a global field as intersections of the corresponding
subsets of local fields. In this case we say that these (global) sets
\emph{satisfy a local-global principle}. Local-global problems are often formulated
in the setting of algebraic varieties. In this case the validity of a local-global
principle implies a description of some subsets of 
the $K$-rational points of certain varieties in terms of their $K_{\mathfrak{p}}$-rational
points.

Along with the Hasse-Minkowski theorem on quadratic forms in Theorem~\ref{HM}, one of the most famous
local-global principles is the Albert-Brauer-Hasse-Noether (ABHN) theorem for central simple algebras.

\begin{theorem}[{Albert, Brauer, Hasse, Noether, 1932, cf.~\cite[Theorem~8.1.17]{NSW}}]\label{thm:ABHN}
Let $K$ be a global field 
and let $A$ be a central simple algebra over $K$. Then  $A$
splits over $K$ if and only if $A$ splits over $K_{\mathfrak{p}}$, for all places $\mathfrak{p}$ of $K$.
Moreover, we have the following five-term exact sequence:

\begin{equation}\label{eq:sesBrauer}
	0\longrightarrow \mathrm{Br}(K) \longrightarrow \bigoplus_{\mathfrak{p} \in M_K} \mathrm{Br}(K_{\mathfrak{p}})\longrightarrow \mathbb{Q}/\mathbb{Z} \longrightarrow 0.
\end{equation}
\end{theorem}

\begin{remark}\label{rmk: explenation of ABHN}
The exact sequence \eqref{eq:sesBrauer} in the ABHN theorem describes the relation of the Brauer group of a global field $K$ with those of all its completions $K_{\mathfrak{p}}$. More precisely:
\begin{enumerate}
    \item The map 
    \[
        \mathrm{Br}(K) \rightarrow \bigoplus_{\mathfrak{p}\in M_K} \mathrm{Br}(K_{\mathfrak{p}})
    \]
    is the diagonal map. Hence, every element $A\in \mathrm{Br}(K)$ has trivial image in $\mathrm{Br}(K_{\mathfrak{p}})$ for all but finitely many places $\mathfrak{p}$. The injectivity of the map is telling us that if an element in $\mathrm{Br}(K)$ has trivial image in $\mathrm{Br}(K_{\mathfrak{p}})$ for all places $\mathfrak{p}$ of $K$, then it has to be trivial. For every field $K$ and every central simple algebra $A/K$ 
    the class of $A$ in $\mathrm{Br}(K)$ is trivial if and only if $A$ is split (see Definition~\ref{def: Brauer group}). Thus, a central simple algebra $A$ over a global field $K$ is split if and only if for every place $\mathfrak{p}$ of $K$ the central simple algebra $A\otimes_K K_{\mathfrak{p}}$ is split over $K_{\mathfrak{p}}$.
    \item The map
    \[
       \bigoplus_{\mathfrak{p}\in M_K} \mathrm{Br}(K_\mathfrak{p})\rightarrow \mathbb{Q}/\mathbb{Z}
    \]
    is given by the sum of the invariant maps $\mathrm{inv}_{\mathfrak{p}}:\mathrm{Br}(K_{\mathfrak{p}})\rightarrow \mathbb{Q}/\mathbb{Z}$. For local fields the invariant map is an isomorphism and its construction can be found in \cite[Chapter XVIII, Section $3$]{LocalFieldsSerre}; for non-archimedean real places it is given by the isomorphism $\mathrm{Br}(\mathbb{R}) \cong \mathbb{Z}/2\mathbb{Z}$ and for non-archimedean complex places it is the zero map. The exactness in the center of the short exact sequence of Theorem \ref{thm:ABHN} tells us that an element $(A_{\mathfrak{p}})\in \oplus_{\mathfrak{p}\in M_K} \mathrm{Br}(K_{\mathfrak{p}})$ comes from an element $A\in \mathrm{Br}(K)$ if and only if $\sum_{\mathfrak{p}\in M_K} \mathrm{inv}_{\mathfrak{p}}(A_{\mathfrak{p}})=0$. Hence, Theorem~\ref{thm:ABHN} gives a recipe to detect whether a (global) element in $\mathrm{Br}(K)$ with prescribed (local) behavior at each completion of $K$ exists.
    \item Finally, if we restrict the short exact sequence to the $2$-torsion elements, then the injectivity of the diagonal map from $\mathrm{Br}(K)[2]$ to $\oplus_{\mathfrak{p}} \mathrm{Br}(K_{\mathfrak{p}})[2]$ follows from Theorem~\ref{HM} (Hasse--Minkowski). In fact, it is proven in \cite[Proposition~1.3.2]{GS06} that a quaternion algebra $(a,b)$ splits (i.e. has trivial image in the Brauer group) if and only if the corresponding conic $ax^2+by^2=z^2$ has a non-trivial solution. 
\end{enumerate}
\end{remark}

\section{Survey}\label{section:survey}

In 1949 Robinson applied the Hasse-Minkowski local-global principle for quadratic forms to show that $\mathbb{Z}$ is definable in $\mathbb{Q}$.
In the decades since then, others have applied a variety of local-global principles to obtain model-theoretic definability results, and this area has been especially active in the last fifteen years.
Most recently the local-global principle for central simple algebras has been used in several significant papers.
In this section we give a brief survey of these results.

\subsection{Robinson's results}\label{Robinson}

In 1949, Julia Robinson proved in \cite{robinson_1949} that the integers are a definable subset of $\mathbb{Q}$:

\begin{theorem}[Robinson, 1949] \label{Rob_Thm_3.1} 
For a rational number $n$ to be an integer it is necessary and sufficient that it satisfies the following formula:

\begin{equation}\label{eq:formula}
\begin{split}
    \forall a,b\Big(\Big\{\Big(\exists x,y,z : 2+bz^2 & =x^2+ay^2 \wedge \forall m \Big[\Big(\exists x,y,z : 2+abm^2+bz^2=x^2+ay^2\Big) \\
    & \rightarrow\Big(\exists x,y,z : 2+ab(m+1)^2+bz^2=x^2+ay^2\Big)\Big]\Big\} \\
    & \rightarrow\Big(\exists x,y,z : 2+abn^2+bz^2=x^2+ay^2\Big)\Big).
\end{split}
\end{equation}

\end{theorem}

\noindent To prove Theorem \ref{Rob_Thm_3.1}, Robinson uses the Hasse-Minkowski theorem \ref{HM} to show that if~$n$ is a rational number such that some of the quadratic forms in \eqref{eq:formula}, with $m=n$, have a solution in $\mathbb{Q}$, then the denominator of $n$ is odd and coprime to various classes of infinitely many prime numbers. She then deduces that the denominator of $n$ is not divisible by any prime number, i.e.\ that $n$ is an integer.

As described in section~\ref{theories},
it follows from Theorem~\ref{Rob_Thm_3.1}
and Gödel's First Incompleteness theorem that the ${\mathfrak{L}}_{\mathrm{ring}}$ theory of $\mathbb{Q}$ is undecidable.  
See for example \cite[Theorems I-IV]{robinson_1949}.

Theorem \ref{Rob_Thm_3.1} also implies the following statement. 

\begin{corollary}[{cf.~\cite[Theorem 3.2]{robinson_1949}}] \label{cor_rob}
Let $n$ be a positive integer and let $R\subseteq\mathbb{Q}^{n}$ be an $n$-ary relation
between rational numbers.
Let $R'\subseteq\mathbb{Z}^{2n}$ be the relation between $2n$ integers
that holds for
$(a_1,\ldots,a_n,b_1,\ldots,b_n)$ if and only if $b_i\neq 0$ for all $i=1,... ,n$ and $R$ holds
for $\left(\tfrac{a_{1}}{b_{n}},\ldots,\tfrac{a_{n}}{b_{n}}\right)$.
Then $R$ is definable in $\mathbb{Q}$ if and only if $R'$ is definable in $\mathbb{Z}$. 
\end{corollary}

\noindent This corollary provides another sense in which $\mathbb{Q}$ and $\mathbb{Z}$ are interpretable in one another:
the map $R\rightarrow R'$ is a bijection between definable subsets of $\mathbb{Q}^{n}$ and certain definable subsets of $\mathbb{Z}^{2n}$.

\subsection{The sets `$S$'}\label{ssec:S}
If $\mathfrak{p}$ is non-archimedean, then $R_{\mathfrak{p}}$ denotes the valuation ring in $K_{\mathfrak{p}}$; we also make use of  the ring $R_{(\mathfrak{p})}:=R_{\mathfrak{p}} \cap K.$ \par
Over a global field $K$ (with $\mathrm{char} \ K \neq 2$) with non-archimedean places $\mathfrak{p}$ and $\mathfrak{q}$, let $a,b,p,q \in K^*$ be such that the quaternion $K$-algebra $H_{a,b}$ ramifies exactly at $\mathfrak{p}$ and $\mathfrak{q}$ and such that 
$\mathrm{ord}_{\mathfrak{p}}p=\mathrm{ord}_{\mathfrak{q}}q=1$ and $\mathrm{ord}_{\mathfrak{p}}q=\mathrm{ord}_{\mathfrak{q}}p=0$.
Eisentr{\"a}ger \cite{Eis05} defines the set (denoted $T$ in \emph{loc.~cit.})
\begin{equation}\label{eq:SabpqE}
\begin{split}
S_{a,b,p,q}^E(K)  :&= \big\{ w\in K \;\big|\;\exists x,y,z\in K \;:\;w^{2}-ax^{2}-by^{2}+abz^{2}= pq \big\} \\
&= \{\; \mathrm{tr}(\alpha)/2 \;\big|\; \alpha \in H_{a,b}(K) : \mathrm{nr}(\alpha) = pq \; \}, 
\end{split}
\end{equation}
i.e. the (halves of) reduced traces of elements of $H_{a,b}$ of reduced norm $pq$.
This set plays an important role in one of her main results, \cite[Theorem~3.1]{Eis05}, which states that for any non-archimedean place~$\mathfrak{p}$ of $K$, the set $R_{(\mathfrak{p})} := \{ x \in K\; \big|\ \mathrm{ord}_{\mathfrak{p}}(x) \geq 0 \}$ (denoted by $R_{\mathfrak{p}}$ in \cite{Eis05}) is diophantine over $K$. 

\begin{remark}
Existential definability of valuation rings of global function fields was proven by Rumely \cite{Rum} and extended to many non-global function fields by Shlapentokh~\cite{Shla_2000}, see also \cite{ShlBook} for other diophantine definitions of valuation rings.
Very recently, Daans showed in his PhD thesis \cite{Daans22} that certain valuation rings of the function field $F$ of a curve are existentially definable, by using local-global principles for quadratic forms
defined over $F$, and that these definitions can be given in a uniform way.
\end{remark}

Following Poonen \cite{Poonen}, for any extension $L/\mathbb{Q}$ 
and $a,b\in L^{\times}$, we define
\begin{equation}\label{eq:SabP}
\begin{split}
S^P_{a,b}(L):&=\big\{2w\in L\;\big|\;\exists x,y,z\in L\;:\;w^{2}-ax^{2}-by^{2}+abz^{2}=1\big\} \\
&= \{\; \mathrm{tr}(\alpha) \;\big|\; \alpha \in H_{a,b}(L) : \mathrm{nr}(\alpha) = 1 \; \}, 
\end{split}
\end{equation}
i.e.\ the set of traces of elements of $H_{a,b}(L)$ of reduced norm $1$.
In fact, Poonen writes $S_{a,b}$ for our $S_{a,b}^{P}(\mathbb{Q})$. He considers the cases $L = \mathbb{Q}$ and $L=\mathbb{Q}_{p}$, and describes the sets $S^{P}_{a,b}(\mathbb{Q}_p)$ in \cite[Lemma~2.1]{Poonen}; importantly, they depend on whether $p \in \Delta_{a,b}:=\Delta_{H_{a,b}}$, i.e., on whether $H_{a,b}$ is split at~$p$.
Precisely, in \cite[Lemma 2.1]{Poonen}, he proves for all rational primes $p$:
\begin{itemize}
	\item
	If $p\not\in\Delta_{a,b}$, then $S^{P}_{a,b}(\mathbb{Q}_{p})=\mathbb{Q}_{p}$;
	\item
	If $p\in\Delta_{a,b}$, then $S^{P}_{a,b}(\mathbb{Q}_{p})\subseteq\mathbb{Z}_{p}$.
 \end{itemize}
 
\begin{remark}\label{S:po} Moreover, the same lemma proves that if $p\in\Delta_{a,b}$ then $S^{P}_{a,b}(\mathbb{Q}_{p})$ is nonempty, because of the containment $\res^{-1}(U_{p}) \subseteq S^{P}_{a,b}(\mathbb{Q}_{p})$, where $\res: \mathbb{Z}_p \to \mathbb{F}_p$ is the reduction-modulo-$p$ map and $U_{p}=\{ \mathrm{Tr}_{\mathbb{F}_{p^2}/\mathbb{F}_p}(\beta) : \mathrm{N}_{\mathbb{F}_{p^2}/\mathbb{F}_p}(\beta) = 1\} \subseteq \mathbb{F}_p$ is the set of finite field traces of elements of norm $1$. 
 \end{remark}

\noindent In addition, he shows (cf.~\cite[Lemma 2.2]{Poonen}) that the sets $S_{a,b}^P$ satisfy a local-global principle if either $a>0$ or $b > 0$:
\begin{equation}\label{eq:LGPforSP}
    S_{a,b}^P(\mathbb{Q}) = \left( \bigcap_{p} S_{a,b}^P(\mathbb{Q}_p) \right) \cap \mathbb{Q},
\end{equation}
where the first intersection is taken over all prime numbers $p$ in $\mathbb{Q}$.\\

Koenigsmann \cite{Koe16} uses the same sets $S_{a,b}^P$ as Poonen over $\mathbb{Q}$ and $\mathbb{Q}_p$. So does Park \cite{Park}, for any number field $K$; in particular, her Lemma~2.2 is the number field analogue of \cite[Lemmas 2.1 and 2.2]{Poonen}. \\

Dittmann \cite{Dit18} extends Poonen's definition of the set $S^P_{a,b}$ from quaternion algebras to any central simple algebra $A$ over any field $K$. More precisely, he defines the set (denoted~$S(A/L)$ in \cite{Dit18})
\begin{equation}\label{eq:SD}
	S^D_{A}(L):=\big\{\mathrm{tr}(x)\;\big|\;x\in A\otimes_K L,  \mathrm{nr}(x)=1\big\},
\end{equation}
for any field extension $L/K$.

In \cite[Lemma~2.5 and Proposition~2.6]{Dit18}, the author describes the sets $S^{D}_{A}$ over local fields. Namely, if $A$ is a central simple algebra of prime degree $l$ over a local field $F$ with valuation ring $\mathcal{O}_F$, then:
\begin{itemize}
	\item
	If $A$ splits over $F$, then $S^{D}_{A}(F)=F$;
	\item
	If $A$ is non-split over $F$ and $l \neq 2$, then $S^{D}_{A}(F)=\mathcal{O}_{F}$.
\end{itemize}
A description is also given in the case $l=2$.
If $K$ is a global field and $A$ is of prime degree, then \cite[Proposition~2.1]{Dit18} gives the local-global principle analogous  to \cite[Lemma 2.2]{Poonen} for the sets $S^D_A$. That is:
\begin{equation}\label{eq:LGPforSD}
    S_{A}^{D}(K) =\left( \bigcap_{\mathfrak{p}} S_{A}^{D}(K_{\mathfrak{p}})\right)\cap K,
\end{equation}
where $\mathfrak{p}$ runs over all places of $K$.\\

Finally, Daans gives a similar, but slightly different, definition of $S$-sets in his \cite{Daans21}; see also \cite{Daans23}. For a field $K$ and a quaternion algebra $Q$ over $K$, he defines 
\begin{equation} \label{Daans:Ssets}
    S_Q^{Da}(K):=\{\mathrm{tr}(x) \ | \ x \in Q\backslash K, \mathrm{nr} (x)=1\}.
\end{equation}
See Subsection \ref{subsec:comparisonS} for a comparison with Poonen's $S_{a,b}^P$.

\subsection{The sets `$T$'} \label{setsT}

In each of the papers surveyed here, the authors define from the sets $S$ a new family of sets $T$.
Importantly, the new sets $T$ should have the same `logical complexity' as the sets $S$:
more precisely, the sets $S$ and $T$ should be definable by existential formulas in the language of rings.\\ 

Following Poonen \cite{Poonen}, for any field extension $L/\mathbb{Q}$, we define the set (denoted $T_{a,b}$ and only considered for $L = \mathbb{Q}$ in \emph{loc.~cit.}) 
\begin{equation}\label{eq:TabP}
T^P_{a,b}(L):=S^P_{a,b}(L)+S^P_{a,b}(L)+\{\;0,1,2,\dots,2309\;\},
\end{equation}
where $a, b \in \mathbb{Q}^{\times}.$ 
In \cite[Lemma~2.5]{Poonen} he proves that for any $a,b\in \mathbb{Q}^\times$ such that 
$a>0$ or $b>0$,
\begin{equation}\label{eq:TP}
    T^P_{a,b}(\mathbb{Q})=\bigcap_{p\in \Delta_{a,b}} \mathbb{Z}_{(p)},
\end{equation}
where $\mathbb{Z}_{(p)}:=\mathbb{Q}\cap \mathbb{Z}_p$, and $\Delta_{a,b}=\Delta_{H_{a,b}}$ as before. 

\begin{remark} \label{2309:po}
For any prime $p>11$, ~\cite[Lemma~2.3]{Poonen} shows that $U_p + U_p = \mathbb{F}_p$ with $U_p$ as in Remark~\ref{S:po}; i.e.~the finite field traces cover all of $\mathbb{F}_p$ for $p$ large enough. To deal with the remaining primes, the set $\{0, 1, 2 \ldots, 2309 = 2\cdot 3 \cdot 5 \cdot 7 \cdot 11 - 1 \}$ appears in~\eqref{eq:TabP}, because for any $t \in \bigcap_{p \in \Delta_{a,b}} \mathbb{Z}_{(p)}$, there exists $n \in \{1,2,\dots, 2309\},$ so that $\mathrm{red}_p(t-n) \in U_p + U_p$ for all $p \leq 11$, as well. This is used to prove the right-to-left inclusion in~\eqref{eq:TP}.  The left-to-right inclusion is a consequence of (\ref{eq:LGPforSP}).
\end{remark}

Regarding $T$-sets, in the spirit of Koenigsmann, for every field extension $L$ of $\mathbb{Q}$ we define
\begin{equation} \label{int:ko}
T^K_{a,b}(L):=S^P_{a,b}(L)+S^P_{a,b}(L).
\end{equation}
In \cite[Definition~4]{Koe16}, Koenigsmann defines this set (denoted $T_{a,b}$ in \emph{loc.~cit.}) only over the field of rational numbers.
Then \cite[Proposition~6]{Koe16} strengthens \cite[Lemma 2.5]{Poonen} by showing that for \emph{any} $a, b \in \mathbb{Q}^\times$,
\begin{equation}\label{eq:TK}
    T^K_{a,b}(\mathbb{Q})=\bigcap_{p\in \Delta_{a,b}} \mathbb{Z}_{(p)}
\end{equation}
where $\mathbb{Z}_{(\infty)}:=\big\{x\in \mathbb{Q} \;\big|\; -4 \leq x\leq 4 \big\}$.
Whenever at least one of $a, b$ are positive, the $T$-set defined in \eqref{eq:TK} coincides with the $T$-set defined by Poonen.
For Poonen, the set $\mathbb{Z}_{(\infty)}$ does not appear, seeing as $\infty \not\in \Delta_{a,b}$ when $a>0$ or $b>0$. See Subsection~\ref{remarkComparingTs} for more details. \\ 

Park~\cite{Park} generalizes the above to number fields $K$ and defines the set $T$ analogously to Koeningsmann: 
\begin{equation}\label{eq:TabPa}
T_{a,b}^{Pa}(K) := S^P_{a,b}(K) + S^P_{a,b}(K).   
\end{equation}
Her \cite[Proposition 2.3]{Park} states that
\begin{equation} \label{lo:park}
    T_{a,b}^{Pa}(K) = \bigcap_{\mathfrak{p} \in \Delta_{a,b}} R_{(\mathfrak{p})}
\end{equation}
for any $a,b\in K^\times$ such that $\sigma(a)>0$ or $\sigma(b)>0$ for any real embedding $\sigma$ of~$K$, where $R_{(\mathfrak{p})} = \{ x \in K : \mathrm{ord}_{\mathfrak{p}}(x) \geq 0 \}$ is denoted by $R_{\mathfrak{p}}$ in \cite{Park} for any non-Archimedean place $\mathfrak{p}$ of $K$.\\

In \cite{Dit18}, Dittmann gives a definition of the set $T$ that depends on the degree of the central simple algebra $A/K$. Let $K$ be a global field, and $A/K$ a central simple algebra of prime degree $l$. Then 
\begin{equation}\label{eq:TDdef}
T^D_A(K):=\begin{cases} S^D_A(K)-S^D_A(K) & \text{ if }l=2;\\
S^D_A(K) & \text{ otherwise.}
\end{cases}
\end{equation}
Then \cite[Proposition~2.9]{Dit18} is the analogue of \cite[Lemma~2.4]{Poonen}; it proves that 
if $A$ splits at all real places of $K$, then
\begin{equation}\label{eq:TD}
    T^D_A(K)=\bigcap_{\mathfrak{p} \in \Delta_{A/K}}R_{(\mathfrak{p})},
\end{equation}
where $\Delta_{A/K}$ is the finite set of places of $K$ at which $A$ does not split. In \cite{ADF2}, the authors define the $T$-sets as in~\eqref{eq:TDdef}.\\

In \cite{Daans21}, the author defines $T$-sets for a quaternion algebra $Q$ over an arbitrary field~$K$ as follows:
\begin{equation} \label{Daa:Tsets}
    T^{Da}_Q(K):=\bigcap_{\mathfrak{p} \in \Delta_Q}R_{(\mathfrak{p})}, 
\end{equation}
where $\Delta_Q$ is the set of discrete valuations on $K$ for which $Q_\mathfrak{p}:=Q \otimes_K K_\mathfrak{p}$ is non-split.
Daans then shows in \cite[Theorem 5.1]{Daans21}, that for a \emph{non-real} quaternion algebra over a global field $K$ (i.e.~a quaternion algebra that splits over all the archimedean places of~$K$) one obtains 
\begin{equation} \label{Daa:Tsets2}
    T^{Da}_Q(K)=S_Q^{Da}(K) + S_Q^{Da}(K).
\end{equation}

\subsection{Constructing diophantine sets using local-global principles}\label{ssec:intersections}

We present an outline of how the Hasse-Minkowski (Theorem~\ref{HM}) and ABHN (Theorem~\ref{thm:ABHN}) local-global principles appear in the papers surveyed here, 
beginning with \cite{Eis05}.
We suppose throughout this section that the characteristic of $K$ is not $2$.\\

\begin{remark}
All the sets $S$ and $T$ appearing in the preceding sections are diophantine subsets over the relevant fields.
For example, $T^{Pa}_{a,b}(K)$ is diophantine over $K$, for each number field $K$.
\end{remark}

\noindent \textbf{LGPs and Eisenträger \cite{Eis05}.} In \cite{Eis05}, Eisenträger proves that given a global field~$K$ and a non-archimedean place $\mathfrak{p}$, the ring 
$R_{(\mathfrak{p})}:=\{x \in K: \mathrm{ord}_{\mathfrak{p}}(x) \geqslant 0\}$ is diophantine over $K$. 
Since
$R_{(\mathfrak{p})}=R_{(\mathfrak{p})} \cap R_{(\mathfrak{q})} + R_{(\mathfrak{p})} \cap R_{(\mathfrak{l})}$ for any two additional non-archimedean places $\mathfrak{q}$ and $\mathfrak{l}$ of $K$,
it suffices to show that $R_{(\mathfrak{p})} \cap R_{(\mathfrak{q})}$ is diophantine over $K$ for every two non-archimedean places~$\mathfrak{p}$ and~$\mathfrak{q}$. This is done in Theorem 3.1 of {\em loc.~cit.},  
 as follows:

\begin{enumerate}
    \item By ABHN, there exists a quaternion $K$-algebra $H_{a,b}$ ramified \emph{precisely} at $\mathfrak{p}$ and $\mathfrak{q}$, see (b) in Remark~\ref{rmk: explenation of ABHN};
    we may assume that $a, b \in R_{(\mathfrak{p})} \cap R_{(\mathfrak{q})}$. As just before Equation~\eqref{eq:SabpqE} defining $S^E_{a,b,p,q}$,  
    let $p,q \in K$ be such that $\mathrm{ord}_{\mathfrak{p}}p=\mathrm{ord}_{\mathfrak{q}}q=1$ and $\mathrm{ord}_{\mathfrak{p}}q=\mathrm{ord}_{\mathfrak{q}}p=0.$ 
    For large enough $r$, Eisentr{\"a}ger proves that $(pq)^r S_{a,b,p,q}^E(K) \subseteq R_{(\mathfrak{p})} \cap R_{(\mathfrak{q})}$.

    \item Let $x \in pR_{(\mathfrak{p})} \cap q R_{(\mathfrak{q})}$ and let $\beta$ be a root of the polynomial $X^2-2xX+pq$. 
    The algebras $H_{a,b}(K_{\mathfrak{p}})$ and $H_{a,b}(K_{\mathfrak{q}})$ now split over the fields $K_{\mathfrak{p}}(\beta)$ and $K_{\mathfrak{q}}(\beta)$, respectively; 
    hence $H_{a,b}(K(\beta))$ is split \emph{everywhere} locally. By ABHN, the algebra $H_{a,b}(K(\beta))$ then splits over $K(\beta)$, see (a) in Remark~\ref{rmk: explenation of ABHN}. 
    It follows that $K(\beta) \subseteq H_{a,b}$, so ${\beta \in H_{a,b}}$ and $x \in S_{a,b,p,q}^E$. Thus, $pR_{(\mathfrak{p})} \cap q R_{(\mathfrak{q})} \subseteq S_{a,b, p, q}^E$, implying that $G:=p^{r+1}R_{(\mathfrak{p})} \cap q^{r+1} R_{(\mathfrak{q})} \subseteq (pq)^rS_{a,b,p,q}^E \subseteq R_{(\mathfrak{p})} \cap R_{(\mathfrak{q})}$ for some $r$, see also point (a) above. 
    
    \item Since $G$ has finite  index in $R_{(\mathfrak{p})} \cap R_{(\mathfrak{q})}$, by introducing coset representatives $s_i$, one obtains $R_{(\mathfrak{p})} \cap R_{(\mathfrak{q})}= \bigsqcup_{i=1}^t (s_i+G) = \bigcup_{i=1}^t (s_i+ (pq)^rS_{a,b, p,q}^E).$ 
Since $S_{a,b, p, q}^E$ is diophantine over $K$ and $R_{(\mathfrak{p})} \cap R_{(\mathfrak{q})}=s_1+s_2+\cdots s_t+(pq)^rS_{a,b,p,q}^E$, the set $R_{(\mathfrak{p})} \cap R_{(\mathfrak{q})}$ is diophantine over $K$. 

In a sense, the set $T^E_{a,b,p,q}:=(pq)^rS^E_{a,b,p,q}+\{  s_1,\dots,s_t \}$ can be seen as a precursor of the $T$-sets of the other authors. Indeed, in Eisenträger's construction the elements $a,b,p,q \in K^\times$ are chosen such that $\Delta_{H_{a,b}}=:\Delta_{a,b}=\{ \mathfrak{p},\mathfrak{q} \}$ and 
\[
    T^E_{a,b,p,q}=R_{(\mathfrak{p})} \cap R_{(\mathfrak{q})}=\bigcap_{\mathfrak{l}\in \Delta_{a,b}} R_{(\mathfrak{l})};
\]
cf.~Equations~\eqref{eq:TP}, \eqref{eq:TK}, \eqref{lo:park} and \eqref{eq:TD}. 
\end{enumerate}

\begin{remark} \label{EisvsPo}
As ABHN was used only in the case of quaternions in the above, it could have been replaced by the Hasse-Minkowski theorem, as in the works of Robinson and Poonen. As mentioned in part (c), in Remark~\ref{rmk: explenation of ABHN} it also follows from the Hasse--Minkowski Theorem~\ref{HM} that a quaternion algebra $H_{a,b}$ over a global field $K$ is unramified if and only if it is unramified over all the completions~$K_{\mathfrak{p}}$. Moreover, Robinson and Poonen work over the field of rational numbers, and for every couple of places $\mathfrak{p}, \mathfrak{q}$, it is possible to construct "by hand" a quaternion algebra $H_{a,b}$ that splits exactly at $\mathfrak{p}$ and~$\mathfrak{q}$.
\end{remark}

\noindent \textbf{LGPs and Poonen \cite{Poonen}.} By definition, for $s\in \mathbb{Q}$, membership of $s$ in 
Poonen's~$S_{a,b}^P$ is equivalent to the existence of an element $x$ of the 
quaternion $\mathbb{Q}$-algebra $H_{a,b}$ such that $\mathrm{tr}(x)=s$ and 
$\mathrm{nr}(x)=1$. Equivalently, the characteristic polynomial of $x$ over~$\mathbb{Q}$ is a quadratic polynomial $X^2-sX+1$. In the local-global 
principle~\eqref{eq:LGPforSP} for $S_{a,b}^P$ where either $a >0$ or $b> 0$, the left-to-right inclusion follows from the definition. To obtain the reverse inclusion, let $s \in 
\mathbb{Q} \cap \bigcap_{p} S_{a,b}^{P}(\mathbb{Q}_p)$. Then for any prime $p$, there exists a solution to $X^2-sX+1=0$ in~$H_{a,b}(\mathbb{Q}_p)$. By the Hasse-Minkowski theorem (Theorem~\ref{HM}), this quadratic polynomial also has a solution in~$H_{a,b}$, thus $s \in S_{a,b}(\mathbb{Q})$.

Poonen proves that $\mathbb{Z}$ is universal-existentially definable in $\mathbb{Q}$ by using the following description of $\mathbb{Z}$ in terms of the $T$-sets $T_{a,b}^P$ (which are in turn defined through the $S$-sets $S_{a,b}^P$): \begin{equation} \label{gl:po}
    \bigcap_{a,b \in \mathbb{Q}_{>0}} T_{a,b}^P(\mathbb{Q})=\mathbb{Z}.
\end{equation}
By Equation~\eqref{eq:TP}, see also Remark \ref{2309:po}, it suffices to show that for any rational prime $p$, there exist rational numbers $a>0$, $b>0$ such that the quaternion $\mathbb{Q}$-algebra $H_{a,b}$ is ramified at $p$. This is a direct consequence of ABHN. However, Poonen avoids using ABHN altogether, by constructing (in a simple manner) such a quaternion algebra ``by hand'' in \cite[Lemma 2.7]{Poonen}: if $p=2,$ take $a=b=7$; otherwise, take $a=p$ and an integer $b$ with $\mathrm{red}_p(b) \in \mathbb{F}_p^{\times} \backslash \mathbb{F}_p^{\times 2}$, where as before, $\mathrm{red}_p$ is the reduction-modulo-$p$ map. 

\begin{remark}
It follows from the definition of $S^E_{a,b,p,q}$ in Equation~\eqref{eq:SabpqE} that in general $S^E_{a,b,p,q}(K) \not\subseteq R_{(\mathfrak{p})}$ for $\mathfrak{p} \in \Delta_{a,b}$; this is remedied by choosing $r$ large enough so that $p^rS^E_{a,b,p,q}(K) \subseteq R_{(\mathfrak{p})}$. 
By contrast, by prescribing reduced norm $1$ when defining $S^P_{a,b}$ in Equation~\eqref{eq:SabP}, it follows that $S_{a,b}^P(\mathbb{Q}_p) \subseteq \mathbb{Z}_p$ for all $p \in \Delta_{a,b}$.
The removal of the factor $p^{r}$ is an important step
towards finding a uniform diophantine definition of $\mathbb{Z}_{(p)}$ in $\mathbb{Q}$,
which is important for obtaining the description of $\mathbb{Z}$ in Equation~\eqref{gl:po}, by seeing $T^P_{a,b}(\mathbb{Q})$ as an intersection of all localizations $\mathbb{Z}_{(p)}$ as in Equation~\eqref{eq:TP}.
\end{remark}

\noindent \textbf{LGPs and Koenigsmann \cite{Koe16}.} In Koenigsmann's \cite{Koe16}, the situation with respect to ABHN is the same as in \cite{Poonen}. Given the description of the $T$-sets $T_{a,b}^K$ via  \eqref{eq:TK},  
it is immediate that \eqref{gl:po} remains true when replacing the sets $T_{a,b}^P$ with the sets $T_{a,b}^K$; cf.~\cite[p.~79]{Koe16}. See also Subsection \ref{remarkComparingTs}.\\

\noindent \textbf{LGPs and Park \cite{Park}.} 
The set $S^{Pa}_{a,b}$ satisfies a
local-global principle in the same way as the set $S^{P}_{a,b}$ does.  
In fact, as described above, for a number field $K,$ since $S^{Pa}_{a,b}(K)$ is defined by a quadratic form
$w^{2}-ax^{2}-by^{2}+abz^{2}-1$, the Hasse-Minkowski theorem assures
an inclusion $K\cap \bigcap_{\mathfrak{p}\in \Delta_{a,b}} S_{a,b}^{Pa}(K_{\mathfrak{p}})\subseteq S_{a,b}^{Pa}(K)$. The reverse inclusion follows by definition, so
\begin{equation} \label{Sabpa_lgp}
S_{a,b}^{Pa}(K)=K\cap \bigcap_{\mathfrak{p}\in \Delta_{a,b}} S_{a,b}^{Pa}(K_{\mathfrak{p}}),
\end{equation}
for all number fields $K$.

As above, let $R_{\mathfrak{p}}$ be the ring of integers of $K_\mathfrak{p}$ and let $R$ be the ring
of integers of $K$. As a consequence of
Equation \eqref{Sabpa_lgp}, we have that $T^{Pa}_{a,b}(K)\subseteq \cap_{\mathfrak{p}\in \Delta_{a,b}} R_\mathfrak{p}$. Let $F_\mathfrak{p}$
be the residue field of $R_\mathfrak{p}$. To prove
the other inclusion, Park considers the set $U_\mathfrak{p}$ of elements $s\in F_\mathfrak{p}$
such that $x^2-sx+1$ is irreducible over~$F_\mathfrak{p}$
and its preimage $V_\mathfrak{p}$ in $R_{\mathfrak{p}}$ (for $\mathfrak{p}$ lying above $p\geq 11$; for $p<11$ the definition 
of $V_\mathfrak{p}$ is slightly different, cf.~\cite[Proposition 2.3]{Park}). 
She then shows that $V_\mathfrak{p}+V_\mathfrak{p}=R_{\mathfrak{p}}$ and uses strong approximation and the local-global
principle given in Equation \eqref{Sabpa_lgp} another time to  get the other inclusion.
Thus
\begin{equation} \label{OK}
\bigcap_{a,b} T^{Pa}_{a,b}(K)= \bigcap_{\mathfrak{p} \in M_{K,na}} R_{(\mathfrak{p})} = R,
\end{equation}
where the intersection runs over $a,b \in K^{\times}$ such that $\sigma(a) > 0$ or $\sigma(b) > 0$ for each real archimedean place $\sigma$, and where $M_{K,na}$ is the subset of $M_k$ containing all the non-archimedean places.

Thus, just as Poonen and Koenigsmann do, 
at this point we may already deduce that the ring of integers $R$ of a number field $K$ is universal-existentially definable in~$K$.
The description of $R$
given in Equation \eqref{OK} is  analogous to that of \cite{Poonen}, cf. (\ref{gl:po}).
She then uses class field theory to obtain a first-order universal formula defining~$R$ in~$K$.\\

\noindent \textbf{LGPs and Dittmann \cite{Dit18}.} In \cite[Proposition 2.1]{Dit18}, Dittmann also uses a local--global principle to prove \eqref{eq:LGPforSD}. Let $K$ be a global field and $A$ a central simple algebra over $K$ of prime degree~$l$. As for Poonen, the left-to-right inclusion in \eqref{eq:LGPforSD} is seen directly from the definitions. We illustrate Dittmann's use of ABHN for proving the reverse inclusion. Let $a \in \bigcap_{\mathfrak{p}} S^D_A(K_{\mathfrak{p}}) \cap K.$
\begin{enumerate}
    \item If $A$ is split over $L$, then \eqref{eq:LGPforSD} is obtained quickly as both sides coincide with $K$. So assume it is not and let $\mathfrak{p}_i$, $i=1,2,\dots, n$ be the places of $K$ where $A$ is ramified. Then for any $i$, there exists $x_{i} \in  A \otimes_K K_{\mathfrak{p}_i}$ of reduced trace $a$ and reduced norm $1$ by assumption. If $x_i \in K_{\mathfrak{p}_i}$, then $x_i^l=1$ by the reduced norm hypothesis, so one obtains $x_i=a/l \in K$ if $l\neq \mathrm{char}(K)$, and $x_i=1$ otherwise, from the reduced trace hypothesis.  Since these quantities are both global, the reverse inclusion follows in this case.
\item Assume $x_i \not \in K_{\mathfrak{p}_i}$, so that $K_{\mathfrak{p}_i}(x_{i})/K_{\mathfrak{p}_i}$ is a degree $l$ field extension splitting $A \otimes_{K} K_{\mathfrak{p}_i}$. Moreover, the minimal polynomial $f_i$ of $x_i$ has degree  $(l-1)$-coefficient equal to $a$ and constant term $(-1)^l$ by assumption. By a weak approximation argument (\cite[Lemma 2.4]{Dit18}), one can assume that the field extensions of the $K_{\mathfrak{p}_i}$, obtained by adjoining a root of the respective $f_i$, all coincide with those obtained by adjoining a root of some common degree~$l$ polynomial $f \in K[X]$. 
Then, by construction, the global field extension $K[X]/(f)$ splits $A$, meaning it is contained in $A$, see (a) in Remark~\ref{rmk: explenation of ABHN}. The global element $X + (f) \in K[X]/(f)$ has reduced trace $a$ and reduced norm $1$, so we are done. 
\end{enumerate}

 We remark that there are common points between the above proof and that of Eisenträger. One difference with \cite{Eis05} and \cite{Poonen} is that in \cite{Dit18} the degree of the central simple algebra is a prime number (possibly larger than $2$), meaning the Hasse-Minkowski theorem is not enough, so the more general ABHN is invoked.

Dittmann then also uses ABHN to prove that for any non-archimedean place $\mathfrak{p}$ of a global field~ $K$, the ring $R_{(\mathfrak{p})}$ is diophantine over $K$. This is done similarly as in Eisenträger's \cite{Eis05} and Poonen's \cite{Poonen}, 
except that instead of constructing a suitable global central simple algebra ``by hand'', ABHN is used to find two central simple algebras with prescribed local splitting behavior, see (b) in Remark~\ref{rmk: explenation of ABHN}. More precisely, one can find central simple algebras $A, A'$ over $K$ such that  $T_A^D(K) + T_{A'}^D(K)=R_{(\mathfrak{p})}$, thus proving the result. A predecessor of \cite{Dit18} is the paper \cite{Rum} produced by
Rumely in 1980, in which the techniques are similar with the use of norm forms instead of quadratic forms (which instead were previously used by Robinson). \\

\noindent \textbf{LGPs and Daans \cite{Daans23}.}
In \cite{Daans23}, Daans picks up on the theme of analyzing the exact number of quantifiers used in definitions of rings of integers in global fields,
following for example
\cite{EM,ZS} and the framework of \cite{DDF}.
Using ABHN he shows that for any finite set~$V$ of non-archimedean places of a global field $K$, the ring $\bigcap_{\mathfrak{p}\in V} R_{(\mathfrak{p})}$ has an $\exists_3$-$\mathfrak{L}_{\mathrm{ring}}(K)$-definition in $K$,
i.e., there is an existential $\mathfrak{L}_{\mathrm{ring}}$-formula $\varphi(x,y_{1},\ldots,y_{m})$ with three existential quantifiers and elements $b_{1},\ldots,b_{m}\in K$ such that 
$K\models\varphi(a,b_{1},\ldots,b_{m})$ if and only if $a\in\bigcap_{\mathfrak{p}\in V}R_{(\mathfrak{p})}$.
In particular, in Lemma~4.1 of \emph{loc.~cit.} the injectivity of the first map in Equation~\eqref{eq:sesBrauer} gives a description of $\bigcap_{\mathfrak{p}\in V} R_{(\mathfrak{p})}$ in terms of the splitting behaviour of a particular quaternion algebra.
Then in Proposition~4.2 of \emph{loc.~cit.} he uses Equation~\eqref{eq:sesBrauer} again to prove that such a quaternion algebra always exists.
In addition, in Section~5 of \emph{loc.~cit.} a global quaternion algebra with prescribed splitting is constructed using ABHN, in order to prove that  for any non-empty finite set $V$ of non-archimedean places, 
the intersection $\bigcap_{\mathfrak{p} \in M_K\setminus V} R_{(\mathfrak{p})}$ is $\forall_{10}$-$\mathfrak{L}_{\mathrm{ring}}(K)$-definable in~$K$;
likewise this means that there is a universal $\mathfrak{L}_{\mathrm{ring}}$-formula $\varphi(x,y_{1},\ldots,y_{m})$ with ten universal quantifiers and elements $b_{1},\ldots,b_{m}\in K$ such that
$K\models\varphi(a,b_{1},\ldots,b_{m})$ if and only if $a\in\bigcap_{\mathfrak{p}\in M_{K}\setminus V}R_{(\mathfrak{p})}$.

\section{Local-global methods and definability}\label{sec:LGPdefinability}

In the previous section, we surveyed a number of recent articles that use local-global principles to answer model-theoretic questions. In essence, the surveyed articles introduce similar notions of sets `$S$', consisting of reduced traces of elements of certain simple algebras with prescribed reduced norm, and construct similar sets `$T$' out of these by taking sums or differences. The sets `$T$' turn out to be intersections of certain other sets; using this observation, they can be showed to be diophantine, by using local-global principles to construct global algebras with suitable properties.
In this section, we make some more detailed comparisons between the methods employed in the surveyed articles in Subsection~\ref{ssec:comp},
give a more model-theoretic outlook in Subsection~\ref{ssec:mod},
and sketch two possible applications of this framework using other local-global principles in Subsection~\ref{ssec:LGP}.

\subsection{Comparisons between surveyed articles}\label{ssec:comp}\
With a view towards both establishing a common model-theoretic framework and finding new applications of other local-global principles, we first provide some more detailed  comparisons between the sets introduced in the surveyed articles.

\subsubsection{Comparisons between the sets `$S$'}\label{subsec:comparisonS}
Let us compare the definitions of $S^E_{a,b,p,q}$, $S^P_{a,b}$ and $S^{Da}_Q$ in Equations~\eqref{eq:SabpqE},\eqref{eq:SabP} and~\eqref{Daans:Ssets}, respectively. 
We start by recalling (see Subsection~\ref{ssec:S}) that the definition of the $S$-sets of Park in \cite{Park}, Koenigsmann in \cite{Koe16} and Dittmann in \cite{Dit18} essentially  coincide with that of Poonen in \cite{Poonen}.

\begin{itemize}
\item $\boldsymbol{S_{a,b, p, q}^E}$ {\bfseries versus} $\boldsymbol{S_{a,b}^P}$:  In 
characteristic different from 2
the difference between taking half of the reduced trace $\mathrm{tr}(\alpha)/2$ in Eisenträger's~\eqref{eq:SabpqE} versus taking $\mathrm{tr}(\alpha)$ in Poonen's~\eqref{eq:SabP} in the definitions of the $S$-sets is inconsequential.
In fact, the proofs of \cite{Poonen}
would still go through if
we instead define $S_{a,b}^P(L) = \{ \mathrm{tr}(\alpha)/2 \;\big|\; \alpha \in H_{a,b}(L) : \mathrm{nr}(\alpha) = 1 \}$ and replace $U_{p}=\{ \mathrm{Tr}_{\mathbb{F}_{p^2}/\mathbb{F}_p}(\beta) : \mathrm{N}_{\mathbb{F}_{p^2}/\mathbb{F}_p}(\beta) = 1\} \subseteq \mathbb{F}_p$ with $U_p/2$, 
cf.~Remark~\ref{S:po}.

\item $\boldsymbol{S_{a,b}^P}$ {\bfseries versus} $\boldsymbol{S^{Da}_{H_{a,b}}}$: While Daans gives a slightly different definition of the $S$-sets in \cite{Daans21}, it is comparable to that of Poonen's (\ref{eq:SabP}). For a quaternion algebra $Q$ over a field $K$, if $x \in K$ is such that $\mathrm{nr}_Q(x)=x^2=1$, then $\mathrm{tr}_Q(x)=2x \in \{\pm 2\}$, where $\mathrm{nr}_Q$ and $\mathrm{tr}_Q$ denote the reduced norm and reduced trace in $Q$, respectively.  
By \cite[Proposition 3.4]{Daans23}, $Q$ is split over $K$ if and only if there exists $\alpha \in Q \backslash K$ such that $\mathrm{nr}_Q(\alpha)=1$ and $\mathrm{tr}_Q(\alpha) \in \{\pm 2\}$, meaning:
\begin{enumerate}
    \item if $Q$ is split over $K$, then 
    $S^{Da}_Q(K):=\{\mathrm{tr}_Q(x) \ | \ x \in Q \backslash K, \mathrm{nr}_Q(x)=1\}=\{\mathrm{tr}_Q(x) \ | \ x \in Q, \mathrm{nr}_Q(x)=1\}$,
\item if $Q$ is non-split over $K$, then $S^{Da}_Q(K)=\{\mathrm{tr}_Q(x) \ | \ x \in Q, \mathrm{nr}_Q(x)=1\} \backslash \{\pm 2\}.$
\end{enumerate}
If $K$ is a field extension of $\mathbb{Q}$, then there exist $a, b \in K^{\times}$ such that $Q=H_{a,b}$. So, if $Q$ is split over $K$, then $S^{Da}_Q(K)=S^P_{a,b}(K)$, and if not, then $S^{Da}_Q(K)=S_{a,b}^P(K) \backslash \{\pm 2\}.$ This difference is however minor, since the corresponding $T$-sets will still coincide, see Subsection \ref{remarkComparingTs}.
\end{itemize}

\subsubsection{Comparisons between the sets `$T$'}\label{remarkComparingTs} 
Let us now compare the definitions of $T^P_{a,b}$, $T^K_{a,b}$, $T^{Pa}_{a,b}$, $T^D_A$ and $T^{Da}_Q$ made in Equations~\eqref{eq:TabP}, \eqref{int:ko}, \eqref{eq:TabPa}, \eqref{eq:TDdef}, and \eqref{Daa:Tsets}, respectively.
\begin{itemize}
\item $\boldsymbol{T_{a,b}^P}$ {\bfseries versus} $\boldsymbol{T_{a,b}^K}$: 
For every $a,b\in \mathbb{Q}^\times$, it follows immediately from the definitions in Equations~\eqref{eq:TabP} and~\eqref{int:ko} that $T_{a,b}^K(\mathbb{Q})\subseteq T_{a,b}^P(\mathbb{Q})$.        Moreover, if $a>0$ or $b>0$, then the quaternion $\mathbb{Q}$-algebra $H_{a,b}$ splits over the real numbers; therefore, Equations~\eqref{eq:TP} and~\eqref{eq:TK} tell us that in this case actually $T^K_{a,b}(\mathbb{Q})=T^P_{a,b}(\mathbb{Q})=\bigcap_{p \in \Delta_{a,b}} \mathbb{Z}_{(p)}$.

The reason for the appearance of the set $\{0,1,2,\ldots, 2309\}$ in the definition of $T^P_{a,b}$ in~\eqref{eq:TabP} is explained in Remark~\ref{2309:po}. Even though Koenigsmann removes this in his definition~\eqref{int:ko} of $T$-sets, Equation~\eqref{eq:TK}, an analogue of Poonen's \eqref{eq:TP}, remains true:
{for any $a, b \in \mathbb{Q}^{\times}$,}
\begin{equation}\label{eq:TPTK}
T^K_{a,b}(\mathbb{Q}) = \bigcap_{p \in \Delta_{a,b}} \mathbb{Z}_{(p)}. 
\end{equation}
{The left-to-right inclusion is a consequence of the definition of $T_{a,b}^K$ and the fact that for $p \in \Delta_{a,b}$ we have $S_{a,b}^K(\mathbb{Q}_p) \subseteq \mathbb{Z}_p$ if $p \neq \infty$ and $S_{a,b}(\mathbb{R}) \subseteq \mathbb{R}$ otherwise (cf. \cite[pg. 77]{Koe16}) and \cite[Lemma 2.1]{Poonen})}.
To prove the right-to-left inclusion, for $p>11$ prime, Koenigsmann uses the same technique as Poonen, cf.~Remark~\ref{2309:po}.
To prove it for $p \leq 11$ or ${p=\infty}$,  
he constructs ``by hand'' intermediary sets~$V_p$ 
such that $\emptyset \neq V_p \subseteq \mathbb{Z}_{p}$, where $\mathbb{Z}_{\infty} := [-4,4] \subseteq \mathbb{R}$. The~$V_p$ are moreover constructed to satisfy $\mathrm{res}^{-1}(U_p) \subseteq V_p$ for all $p$ prime, where, as in Remark~\ref{S:po}, $\res: \mathbb{Z}_p \to \mathbb{F}_p$ is the reduction-modulo-$p$ map and $U_{p}=\{ \mathrm{Tr}_{\mathbb{F}_{p^2}/\mathbb{F}_p}(\beta) : \mathrm{N}_{\mathbb{F}_{p^2}/\mathbb{F}_p}(\beta) = 1\} \subseteq \mathbb{F}_p$ is the set of finite field traces of elements of norm~$1$.
Then he checks case-by-case that they also satisfy the property $V_p + V_p = \mathbb{Z}_p$ for all~$p$ (including $\infty$). Since $V_p \subseteq S^P_{a,b}(\mathbb{Q}_p)$ is an open subset, this last property ensures that $\bigcap_{p \in \Delta_{a,b}} \mathbb{Z}_{(p)} \subseteq T^K_{a,b}(\mathbb{Q})$, thus proving~\eqref{eq:TPTK}.

\item $\boldsymbol{T_{a,b}^P}$ {\bfseries versus} $\boldsymbol{T_{a,b}^{Pa}}$:  
Proposition~2.3 in~\cite{Park}, cf.~\eqref{lo:park}, is a generalization to number fields~$K$ of Lemma~2.4 of Poonen, cf.~\eqref{eq:TP}. Poonen's sign condition that either $a> 0$ or $b>0$ assures that the quaternion algebra $H_{a,b}$ splits over the real numbers, i.e., at the only archimedean place of~$\mathbb{Q}$. Park's condition, that $\sigma(a)>0$ or $\sigma(b)>0$ for each real archimedean place of $K$, is taking care of the same phenomenon.
Indeed, her condition ensures that the quaternion algebra $H_{a,b}$ splits over all the real places of $K$; any quaternion algebra splits at all complex places.
{See also Section \ref{ssec:intersections} for an idea of the proof which is similar to the one described in the previous bullet point.} 

\item $\boldsymbol{T_{a,b}^K}$ {\bfseries versus} $\boldsymbol{T_A^D}$:  
To define the set $T_{a,b}^K$, Koeningsmann takes the sum of $S^P_{a,b}$ with itself, while to define $T_A^D$, Dittmann takes the difference of $S^D_A$ with itself. Choosing our central simple algebra $A = H_{a,b}$ to be a quaternion algebra over some extension $K$ of $\mathbb{Q}$, the sets $S^P_{a,b}(K)$ and $S^D_A(K)$ coincide. Then to compare $T^K_{a,b}(K)$ and $T^D_A(K)$, it suffices to check that $S^P_{a,b}(K) = S^D_A(K)$ is invariant under multiplication by $-1$.
Recall that $s\in S^P_{a,b}(K)$ if and only if there exists $\alpha\in H_{a,b}$ such that $\mathrm{tr}(\alpha)=s$ and $\mathrm{nr}(\alpha)=1$.
Moreover, $\mathrm{tr}(-\alpha)=-\mathrm{tr}(\alpha)$ and $\mathrm{nr}(-\alpha)=\mathrm{nr}(\alpha)=1$, hence $-s\in S^P_{a,b}(K)$ as well.  In general, the difference of two sets has a topological interpretation: if a set $V$ has non-empty interior in a topological ring $R$, then $V-V$ is a neighbourhood of zero.

 \item    $\boldsymbol{T_{a,b}^{Pa}}$ {\bfseries versus} $\boldsymbol{T_{H_{a,b}}^{Da}}$: Let $Q$ be a non-real quaternion algebra over a number field~$K$. Then there exist $a,b \in K^{\times}$ such that $Q \cong H_{a,b}$, and for any real archimedean place $\sigma$ of $K$, either $\sigma(a)>0$ or $\sigma(b)>0$. Considering Park's relation (\ref{lo:park}), by Daans' definition (\ref{Daa:Tsets}), $T_{a,b}^{Pa}(K)=T_{Q}^{Da}(K)$.

    \end{itemize}

\subsection{From diophantine sets to diophantine families}\label{ssec:mod} 

Although the sets $S=S_{a,b,p,q}^{E},S_{a,b}^{P}$, etc., and $T=T_{a,b}^{P},T_{a,b}^{K}$, etc.,
are often defined differently in the various papers, as described above,
there are clear commonalities.
Most obviously, they are all 
diophantine sets,
or rather families of diophantine sets defined with parameters describing the relevant quaternion algebras and central simple algebras.
That is, for each of the 
$S$-sets
$S$,
in the papers above
(and similarly for the $T$-sets),
there is a field $K$, an existential $\mathfrak{L}_{\mathrm{ring}}$-formula
$\varphi(x,y_{1},\ldots,y_{n})$,
and there are elements $b_{1},\ldots,b_{n}\in K$ 
such that
\[S(F)=\{a\in F\mid F\models\varphi(a,b_{1},\ldots,b_{n})\},\]
for each $F/K$.
This slightly shifts our perspective from a set that is existentially definable in a given field $K$,
to the family of subsets of field extensions of $K$ that are defined by the same formula, using constants from $K$.
More precisely, we recall the following definition from
\cite[Section 3]{ADF2}:
for a field $K$, a {\em $1$-dimensional diophantine family} over $K$
is a map $D$ from the class of field extensions $F$ of $K$ to sets
given by polynomials $f_{1},\ldots,f_{r}\in K[X,Y_{1},\ldots,Y_{m}]$,
for some $m$,
that is:
\[D(F)=\{x\in F\mid\exists y_{1}\ldots\exists y_{m}\;f_{1}(x,y_{1},\ldots,y_{m})=\ldots=f_{r}(x,y_{1},\ldots,y_{m})=0\},\]
for every $F/K$.
From another point of view, such diophantine families are the subsets of fields $F$ extending $K$ that are defined by a given existential formula
\[\varphi(x)=\exists y_{1}\ldots\exists y_{m}\;f_{1}(x,y_{1},\ldots,y_{m})=\ldots=f_{r}(x,y_{1},\ldots,y_{m})=0.\]
The reader will note a slight divergence from the 
nomenclature of Section~\ref{section:logic} in that we allow $r>1$ in the above definition:
$1$-dimensional diophantine families are really families of existentially definable sets,
however according to Remark~\ref{rem:conjunction_elimination}
this rarely makes a difference. 

In the preceding sections we discussed the various sets $S$ and $T$:
they both form $1$-dimensional diophantine families.
For example, Poonen's sets $S_{a,b}^{P}$ form a $1$-dimensional diophantine family over~$\mathbb{Q}$.
In another case, Dittmann's sets $S_{A}^{D}$ form a $1$-dimensional diophantine family over $K$, where $K$ is the field over which the central simple algebra $A$ is defined.
Moreover, each of the sets $S$ considered in~Subsection~\ref{ssec:S} satisfies a local-global principle, namely that
$S(K)=K\cap\bigcap_{\mathfrak{p}\in\Delta}S(K_{\mathfrak{p}})$,
for a finite set $\Delta$ of places of $K$.
Furthermore, each such $S$ has a very simple description in each closure $K_{\mathfrak{p}}$:
$
S(K_{\mathfrak{p}})=K_{\mathfrak{p}}$
if $\mathfrak{p}\notin\Delta$
and otherwise
$S(K_{\mathfrak{p}})\pm S(K_{\mathfrak{p}})=R_{\mathfrak{p}}$,
with the addition of a finite set in the case of Poonen,
as in \eqref{eq:TabP}.

This motivates the general study of $1$-dimensional diophantine families satisfying similar local-global properties:
Let $D$ be a $1$-dimensional diophantine family over a field $K$,
let $\Delta\subseteq M(K)$ be a finite set of (non-archimedean) places of $K$,
let $\tau=(e,f)\in\mathbb{N}^{2}$,
and let $\mathcal{F}_{K}$ be a family
of field extensions of $K$.
Now, for any field extension $F/K$, not necessarily in $\mathcal{F}_K$, let $M_{\Delta}^{\tau}(F)$ be the set of places
$\mathfrak{P}$ of $F$ lying above some place $\mathfrak{p}$ of $\Delta$
of relative type $\leq\tau$, i.e.,~the ramification degree of the extension $\mathfrak{P}/\mathfrak{p}$ is at most $e$ and the inertia degree divides $f$.
Let $R_{\Delta}^{\tau}(F)$ be the holomorphy ring 
\begin{align*}
	\bigcap_{\mathfrak{P}\in M_{\Delta}^{\tau}(F)}R_{(\mathfrak{P})}
\end{align*}
corresponding to $M_{\Delta}^{\tau}(F)$.
For more background on this setting see for example \cite[Section 2]{ADF2}.
Suppose that
\begin{enumerate}
\item
$D(F)=F\cap\bigcap_{\mathfrak{P}\in M(F)}D(F_{\mathfrak{P}})$,
\item
$D(F_{\mathfrak{P}})=F_{\mathfrak{P}}$ if $\mathfrak{P}\notin M_{\Delta}^{\tau}(F)$, and
\item
$D(F_{\mathfrak{P}})=R_{\mathfrak{P}}$ if $\mathfrak{P}\in M_{\Delta}^{\tau}(F)$,
\end{enumerate}
for all extensions $F\in\mathcal{F}_{K}$.
Then $D(F)=R_{\Delta}^{\tau}(F)$ for $F\in\mathcal{F}_{K}$.
This is the perspective of
\cite{ADF2},
with $K$ a number field,
$\mathcal{F}_{K}$ the finite extensions of $K$,
and $\Delta=\{\mathfrak{p}\}$ is a singleton,
to prove that 
$R_{\mathfrak{p}}^{\tau}(F)$ is a $1$-dimensional diophantine family over $K$,
for finite extensions $F/K$
\cite[Proposition 4.8]{ADF2}.
This approach can be immediately generalized to allow finite sets $\Delta$ of primes on $K$,
or to allow $n$-dimensional diophantine families for $n\geq1$.
In another direction, we may replace $M_{\Delta}^{\tau}(F)$ with the family $M_{\Delta}^{\tau}(F;a)=\{\mathfrak{P}\in M_{\Delta}^{\tau}(F)\mid a\in\mathcal{O}_{\mathfrak{P}}\}$ of open-closed subsets of $M_{\Delta}^{\tau}(F)$, given by a fixed $a\in K$.
In \cite[Theorem 6.5, Corollary 6.6]{ADF2} a sufficient condition is given for the holomorphy rings
associated to $M_{\Delta}^{\tau}(F;a)$ to form a $1$-dimensional diophantine family over $K$.
Finally, turning away from holomorphy rings, we may generalize condition (3), for example by instead requiring that $D(F_{\mathfrak{P}})$ be the subset of $F_{\mathfrak{P}}$ defined by a given formula $\varphi(x)$ in the language of rings.

These considerations prompt the search for a general approach to understand the theory of diophantine families $D$ that satisfy (1), in a suitable generality.
For the sake of concreteness, let us consider the case where $K$ has characteristic zero and $\mathcal{F}_{K}$ is just the singleton $\{K\}$.
For our purposes, a {\em classical place} on $K$ is a place $\mathfrak{p}$ lying over a $p$-adic place on $\mathbb{Q}$, and of relative type $\leq \tau$ for some $\tau\in\mathbb{N}^{2}$.
Let $K^{\mathfrak{p}}$ denote a choice of $p$-adic closure of $K$ with respect to $\mathfrak{p}$ -- this is a substitute for the notion of completion. 
Let $M^{\mathrm{c}}(K)$ be the space of classical places on $K$, equipped with the constructible topology given by subbasic open sets 
$\{\mathfrak{p}\in M^{\mathrm{c}}(K)\mid a\in\mathcal{O}_{\mathfrak{p}}\}$
and 
$\{\mathfrak{p}\in M^{\mathrm{c}}(K)\mid a\notin\mathcal{O}_{\mathfrak{p}}\}$,
for $a\in K$.
Let $\mathcal{D}$ be the set of $1$-dimensional diophantine families $D$ over $K$ for which there exists a compact set $M_{D}\subseteq M^{\mathrm{c}}(K)$
with
$M_{D}$ open-closed in $M^{\mathrm{c}}(K)$,
such that
$D(K)=K\cap\bigcap_{\mathfrak{p}\in M_{D}}D(K^{\mathfrak{p}})$.
Then, instead of the language of rings,
we consider the first-order language $\mathfrak{L}$ equipped with a unary predicate symbol $P_{D}$
for each $D\in\mathcal{D}$.
We then form an $\mathfrak{L}$-structure $\mathcal{K}=(K,(D(K))_{D\in\mathcal{D}})$ with underlying set $K$ in which each $P_{D}$ is interpreted in the obvious way by $D(K)$.
In a natural sense $\mathcal{K}$ is less complicated than the field $K$:
each new predicate in $\mathcal{K}$ is definable by a formula in the language of rings (in fact, by an existential formula).
In model-theoretic language, we say that $\mathcal{K}$ is a {\em reduct} of the field $K$.
We can be more precise:
each positive existential $\mathfrak{L}$-formula defines in $\mathcal{K}$ a set that is already existentially definable in the field $K$,
thus the positive existential theory of $\mathcal{K}$ is a reduct of the existential $\mathfrak{L}_{\mathrm{ring}}$-theory of $K$.

\begin{question} \label{question:1}
Under what conditions on $K$, and relative to what information, can one axiomatize the $\mathfrak{L}$-theory of $\mathcal{K}$?
\end{question}

For number fields this problem seems to be extremely hard,
since even to describe $\mathcal{D}$ amounts to deciding which varieties admit some form of local-global principle.
On the other hand, pseudo-algebraically closed fields, pseudo-real closed fields,
pseudo-$p$-adically closed fields,
and others,
satisfy local-global principles for various classes of varieties:
the complete $\mathfrak{L}_{\mathrm{ring}}$-theories of such fields are already understood and there is a large literature.
For example see \cite[Chapter 12]{FJ}, \cite{Prestel}, and \cite{Fehm}.
In such cases, it seems reasonable to expect that the $\mathfrak{L}$-theory of $\mathcal{K}$ can be axiomatized in a similar way.
Perhaps it is possible to find a field $K$ -- outside of the classes mentioned previously -- for which we can axiomatize $\mathcal{K}$, and in so doing we might find new classes of fields of wider model theoretic interest. 
Further study of these structures $\mathcal{K}$ may well shed light on the possibilities and limitations of the local-global method in the analysis of diophantine sets in global fields.

\subsection{Other local-global principles}\label{ssec:LGP}

As described in Subsection~\ref{ssec:intersections}, two local-global principles were used repeatedly in the surveyed articles to obtain diophantine sets; namely, the Hasse-Minkowski and Albert-Brauer-Hasse-Noether theorems, cf.~Theorem~\ref{HM} and~Theorem~\ref{thm:ABHN}, respectively. In this subsection, we put forward two other local-global principles, both arising from the study of elliptic curves, which can potentially be used in a similar way to obtain other diophantine sets. For a comprehensive introduction to the arithmetic of elliptic curves, cf.~e.g.~\cite{Cassels, AEC}.

Elliptic curves have already been used in the context of Hilbert's Tenth problem for global fields. For example, see \cite{CornelissenZahidi}, where the authors use a conjecture on the properties of elliptic divisibility sequences on an elliptic curve to deduce some results on the complexity of defining $\mathbb{Z}$ in $\mathbb{Q}$; see also \cite{Shla_2008, Shla_2012, MRFB}.

\subsubsection{Endomorphism rings of elliptic curves} Let $K$ be any field and consider elliptic curves over~$K$.

An \emph{isogeny} between elliptic curves $E_1,E_2$ is a $K$-morphism  $E_1 \to E_2$  of algebraic groups.  Any nonzero isogeny is surjective and hence has a finite kernel. The set $\mathrm{Hom}_K(E_1,E_2)$ of isogenies from $E_1$ to $E_2$ inherits a group structure from the group structure on~$E_2$, by pointwise addition: $(\varphi+\psi)(P) := \varphi(P) + \psi(P)$. We will denote the inverse of $\varphi$ in this group by $-\varphi$. There is a natural action of $\mathbb{Z}$ on $\mathrm{End}_K(E_1, E_2)$ given by $(n,\varphi) \mapsto n\cdot \varphi$, where $n \cdot \varphi$ means the sum of $|n|$ copies of $\mathrm{sgn}(n) \varphi$, making $\mathrm{End}_K(E_1, E_2)$ a $\mathbb{Z}$-module. 

For an elliptic curve $E/K$, its \emph{endomorphism ring} $\mathrm{End}_K(E) := \mathrm{Hom}_K(E,E)$ has an additional multiplication structure from composing isogenies.  It has no zero divisors and is a torsion-free $\mathbb{Z}$-module. Moreover, it is a \emph{$\mathbb{Z}$-order} in the \emph{endomorphism algebra} $\mathrm{End}_K^0(E) := \mathrm{End}_K(E) \otimes_{\mathbb{Z}} \mathbb{Q}$, which means it is a lattice (i.e., a finitely generated $\mathbb{Z}$-module whose scalar extension to $\mathbb{Q}$ spans $\mathrm{End}^0_K(E)$) that is also a subring.  As a consequence, one shows that the endomorphism ring $\mathrm{End}_K(E)$ is either~$\mathbb{Z}$, an order in a quadratic imaginary field, or an order in a quaternion algebra. The latter moreover only occurs when $\mathrm{char}(K) >0$ (since $\mathrm{End}^0_K(E)$ is commutative when $\mathrm{char}(K) = 0$); we call the elliptic curve \emph{supersingular} in this case. See Sections 10.1 and 10.2 of \cite{Quats} for an introduction to lattices and orders, respectively.

We may localize any $\mathbb{Z}$-lattice $L$ in a $\mathbb{Q}$-vector space $V$ at any prime $\ell$ by forming $L_{\ell} := L \otimes_{\mathbb{Z}}~\mathbb{Z}_{\ell}$, which is a $\mathbb{Z}_{\ell}$-lattice in the $\mathbb{Q}_{\ell}$-vector space $V \otimes_{\mathbb{Q}} \mathbb{Q}_{\ell}$. For any prime $\ell$, we have a natural injection of~$L$ into $L_{\ell}$.
Lattices over $\mathbb{Z}$ (or any Dedekind domain) satisfy a local-global principle \cite[Theorem~9.4.9 and Lemma~9.5.3]{Quats}: given $L$, the map $N \mapsto (N_{\ell})_{\ell-\mathrm{prime}}$ determines a bijection between $\mathbb{Z}$-lattices $N$ in $V$ and the collection of tuples $(N_{\ell})_{\ell}$ of $\mathbb{Z}_{\ell}$-lattices in $V \otimes_{\mathbb{Q}} \mathbb{Q}_{\ell}$ satisfying $L_{\ell}=N_{\ell}$ for all but finitely many primes $\ell$. We note here that the choice of the lattice $L$ is irrelevant, as for any other such $L'$, one has $L_{\ell}=L'_{\ell}$ for all but finitely many primes $\ell$, so there is no ``real" dependence on $L$; see also the comment and remark after \cite[Theorem~9.4.9]{Quats}.

Similarly, we may localize a $\mathbb{Z}$-order $\mathcal{O}$ in a   $\mathbb{Q}$-algebra $R$ at any prime~$\ell$ by forming $\mathcal{O}_{\ell} := \mathcal{O} \otimes_{\mathbb{Z}} \mathbb{Z}_{\ell}$, which is now a $\mathbb{Z}_{\ell}$-order in the $\mathbb{Q}_{\ell}$-algebra $R \otimes_{\mathbb{Q}} \mathbb{Q}_{\ell}$; again, for any prime $\ell$ the order $\mathcal{O}$ naturally embeds into $\mathcal{O}_\ell$. Moreover, for $\mathbb{Z}$-lattices, being an order is a local property, i.e., checking that a given lattice is an order can be done locally at each prime. Hence, the lobal-global principle for $\mathbb{Z}$-lattices above specializes to a local-global principle for $\mathbb{Z}$-orders in~$R$.

For $\ell \neq \mathrm{char}(K)$ a prime, the \emph{$\ell$-adic Tate module} (or \emph{$\ell$-divisible group}) of an elliptic curve $E/K$ is $T_{\ell}(E) := \varprojlim_{n} E[\ell^n]$, where the inverse limit is taken with respect to the natural multiplication-by-$\ell$ maps. The Galois group $\mathrm{Gal}(\bar{K}/K)$ acts on each $E[\ell^n]$ in a way which respects the inverse limit, and hence extends to the Tate module $T_{\ell}(E)$; cf. \cite[III.7]{AEC}. Due to Tate's theorem, localizations of endomorphism rings, viewed as orders in their respective endomorphism algebras, have an interpretation in terms of Tate modules. More precisely, for any two elliptic curves $E_1,E_2$, there is a natural injective map $\mathrm{Hom}_K(E_1,E_2) \otimes_{\mathbb{Z}} \mathbb{Z}_{\ell} \hookrightarrow \mathrm{Hom}(T_{\ell}(E_1),T_{\ell}(E_2))$, where the right hand side consists of  $\mathrm{Gal}(\bar{K}/K)$-\emph{equivariant} morphisms of $\mathbb{Z}_{\ell}$-modules (cf. Definition 34.1.60 of \cite{Quats}). Tate's theorem shows that this map is an isomorphism when $K$ is a finite field; Faltings showed the same when $K$ is a number field. 

When $\mathrm{char}(K) =: p > 0$ and $\ell = p$, instead of the Tate module of $E$ one takes its \emph{(contravariant) Dieudonné module} $M(E)$, or (categorically) equivalently, its associated \emph{$p$-divisible group}; cf.~\cite[\S 1.2]{Wat} for both of these notions. 
When $K$ is moreover finite, Tate also showed that $\mathrm{Hom}_K(E_1,E_2) \otimes_{\mathbb{Z}} \mathbb{Z}_p$ is isomorphic to $\mathrm{Hom}(M(E_2), M(E_1))$.\\

Putting the above together, for an elliptic curve $E$ over a finite or number field $K$ we have an exact sequence
\begin{equation}\label{eq:LGPend}
0 \to \mathrm{End}_K(E) \to \prod_{\ell} \mathrm{End}_K(E) \otimes \mathbb{Z}_{\ell} \simeq \prod_{\ell} \mathrm{End}(T_{\ell}(E)),
\end{equation}
where the product runs over all primes $\ell$, including possibly $\ell = \mathrm{char}(K)$, though in that case $T_{\ell}(E):=M(E)$ is the Dieudonné module. The image of the second map lands in the restricted product where all but finitely many entries are the identity map. 
\\

Inspired by the use of local-global principles for quaternion algebras in the surveyed articles, we aim to use the local-global behavior of endomorphism rings of elliptic curves over finite fields and number fields mentioned above to obtain diophantine sets. To do this, it is necessary to first determine the next term in the exact sequence~\eqref{eq:LGPend}, since the above-mentioned restricted product is hard to work with. It should be noted that while endomorphism rings are orders that have reduced trace and norm maps, the exact sequence~\eqref{eq:LGPend} is a map of \emph{elements} of orders, so a suitable replacement for trace and norm maps is needed in order to define the analogues of the sets~$S$ and~$T$ from the articles we surveyed in Section~\ref{section:survey}.

\subsubsection{Torsion points on elliptic curves}
In this subsection, we let $K$ be a number field. As in the previous sections, we denote by $M_K$ the set of all places of $K$.

Let $A$  be a commutative connected algebraic group defined over $K$. In \cite{DZ}, the authors propose a local-to-global approach to study the divisibility of $K$-points of $A$ by integers. More precisely, they ask the following:

\begin{question} 
Let $m \in \mathbb{N}$ and $P \in A(K)$. Suppose that for all but finitely many $\mathfrak{p} \in M_K$, there exists $Q_\mathfrak{p} \in A(K_\mathfrak{p})$ such that $P=mQ_\mathfrak{p}$ in $A(K_\mathfrak{p})$. Does this imply there exists $Q \in A(K)$ such that $P=mQ$ in $A(K)$?
\end{question}

If the answer to this question is positive, we say that $A$ satisfies a local-global principle for $m$-divisibility.
One can quickly reduce to the case where $m=p^n$ for $p$ a prime number and $n \in \mathbb{N}$. There have been many recent advances in the study of this question; for more details the interested reader can consult the recent survey \cite{DP}.  One of the main tools is to translate the question into a group cohomology problem, where the
obstruction to the local-global principle for $m$-divisibility is given by a subgroup of $H^1(\mathrm{Gal}(K(A[m])/K), A[m])$, denoted by $H^1_{\mathrm{loc}}(\mathrm{Gal}(K(A[m])/K), A[m])$, considered in \cite{DZ}. Similar groups
were introduced before by Tate (as stated by Serre in \cite{Ser_64}) and Sansuc \cite{Sansuc}. The group $H^1(\mathrm{Gal}(K(A[m])/K), A[m])$  contains a subgroup which is isomorphic to the Tate-Shafarevich
group $\Sha(K, A[p^n])$ (see \cite{DP} for further details).
The vanishing of $H^1_{\mathrm{loc}}(\mathrm{Gal}(K(A[m])/K), A[m])$ provides a sufficient condition for the local-global principle to hold \cite[Proposition 2.1]{DZ}. There are also known examples where it does not hold \cite{DZ_04, DZ2, Pal3}. 

In the case where $A=E$ is an elliptic curve, it was shown in \cite[Theorem~3.1]{DZ} that the local-global principle for $m$-divisibility holds for $m=p$ a prime number. 
Later in \cite{PRV_12}, the authors proved this for all powers $m=p^n$ of every prime number $p>\left(3^{\frac{[K:\mathbb{Q}]}{2}}+1\right)^2$.
The bound has been extended to $p\geq 5$  when $K=\mathbb{Q}$ in \cite{PRV_14} and this is best possible because of counterexamples produced for powers of $p=2,3$ in \cite{Pal3, Creutz}. 
In addition, in the special case when~$P$ is a torsion point of an elliptic curve $E$ defined over a number field $K$,
it was shown in \cite{GR} that the local-global principle for $p^n$-divisibility holds for all $p\geq 3$ and $n\geq 1$.

By the definition of the Tate-Shafarevich group, 
one can write down the following exact sequence: 

\begin{equation} \label{eq:Sha}
0 \rightarrow \Sha(K, E[p^n]) \rightarrow H^1(K, E[p^n]) \xrightarrow[]{\pi} \prod_{\mathfrak{p} \in M_K} H^1(K_\mathfrak{p}, E[p^n]),
\end{equation}

which holds whenever $E$ is a commutative algebraic group and in particular
when $E$ is an elliptic curve. 

 Taking the articles surveyed in Section~\ref{section:survey} as motivation, we propose to use the local-global principle for $p^n$-divisibility, through the exact sequence (\ref{eq:Sha}), for the purpose of proving that certain sets, naturally occurring in arithmetic, are diophantine.
In particular, we would like to work in a setting where the local-global principles holds, so that $\Sha(K,E[p^n])$ vanishes. A first crucial step would then be to understand the image of the map $\pi$ in~\eqref{eq:Sha}, which would allow us to determine its cokernel and hence complete the sequence to a short exact sequence. It seems however that, so far, the cokernel of $\pi$ is not well understood.

Another important point is to choose a suitable interpretation of $H^1(K,E[p^n])$ and hence suitable presentations of its elements, e.g. as $p^n$-coverings of~$E$, or $E[n]$-torsors, or theta groups (cf.~\cite[\S 1]{CFONSS}). This choice will determine the context in which the analogues of the sets~$S$ and~$T$ need to be defined and hence in which way the model-theoretic arguments can be adapted and extended. \\

\vskip 0.5cm
\centerline{\textsc{Acknowledgements}}

\par\vskip 0.1cm\noindent This collaboration started at the workshop Women in Numbers Europe 4 which took place in Utrecht, the Netherlands in August 2022. S.~A.\ was supported by GeoMod AAPG2019 (ANR-DFG). V.~K.\ was supported by the Dutch Research Council (NWO) through grant VI.Veni.192.038. Z. K.\ was partially supported by the research training group "GRK 2240: Algebro-Geometric Methods in Algebra, Arithmetic and Topology", funded by the DFG.  L. P.\ is a member of INdAM-GNSAGA.  We thank the referees for helpful comments.


\begin{thebibliography}{10}

\bibitem{ADF2}
Sylvy Anscombe, Philip Dittmann, and Arno Fehm, \emph{A {$p$}-adic analogue of
  {S}iegel's theorem on sums of squares}, Math. Nachr. \textbf{293} (2020),
  no.~8, pp.~1434--1451.

\bibitem{AxKochen}
James Ax and Simon Kochen, \emph{Diophantine problems over local fields. {II}.
  {A} complete set of axioms for {$p$}-adic number theory}, Amer. J. Math.
  \textbf{87} (1965), pp.~631--648.

\bibitem{Bar}
Stefan Bara\'{n}czuk, \emph{On a dynamical local-global principle in
  {M}ordell-{W}eil type groups}, Expo. Math. \textbf{35} (2017), no.~2,
  pp.~206--211.

\bibitem{EGA}
Nicolas Bourbaki, \emph{El\'{e}ments de math\'{e}matique. {XIII}. {P}remi\`ere
  partie: {L}es structures fondamentales de l'analyse. {L}ivre {II}:
  {A}lg\`ebre, {M}odules et {A}nneax semi-simples. {C}hapitre {VIII}},
  Actualit\'{e}s Scientifiques et Industrielles, no. 1261, Hermann \& Cie,
  Paris, 1958.

\bibitem{Cassels}
John W.~S. Cassels, \emph{Lectures on elliptic curves}, London Mathematical
  Society Student Texts, vol.~24, Cambridge University Press, Cambridge, 1991.
  
\bibitem{CHHKPS}
Jean-Louis Colliot-Th\'{e}l\`{e}ne, David Harbater, Julia Hartmann, Daniel Krashen,
  Raman Parimala, and Venapally Suresh, \emph{Local-global principles for zero-cycles on
  homogeneous spaces over arithmetic function fields}, Trans. Amer. Math. Soc.
  \textbf{372} (2019), no.~8, pp.~5263--5286.

\bibitem{CHHKPS_2}
\bysame, \emph{Local-global principles
  for tori over arithmetic curves}, Algebr. Geom. \textbf{7} (2020), no.~5,
  pp.~607--633.

  \bibitem{CHHKPS_3}
\bysame, \emph{Local-global principles for
  constant reductive groups over semi-global fields}, Michigan Math. J.
  \textbf{72} (2022), pp.~77--144.

\bibitem{CoSko}
Jean-Louis Colliot-Th\'{e}l\`ene and Alexei~N. Skorobogatov, \emph{The
  {B}rauer-{G}rothendieck group}, Ergebnisse der Mathematik und ihrer
  Grenzgebiete. 3. Folge. A Series of Modern Surveys in Mathematics, vol.~71,
  Springer, Cham, 2021.

\bibitem{CornelissenZahidi}
Gunther {C}ornelissen and {K}arim {Z}ahidi, \emph{Elliptic divisibility sequences and undecidable problems about rational points}, J. Reine Angew. Math. \textbf{613} (2007), pp.~1--33.

\bibitem{CFONSS}
John~E. Cremona, Tom~A. Fisher, Catherine O'Neil, Denis Simon, and Michael
  Stoll, \emph{Explicit {$n$}-descent on elliptic curves. {I}. {A}lgebra}, J.
  Reine Angew. Math. \textbf{615} (2008), pp.~121--155.

\bibitem{Creutz}
Brendan Creutz, \emph{On the local-global principle for divisibility in the
  cohomology of elliptic curves}, Math. Res. Lett. \textbf{23} (2016), no.~2,
  pp.~377--387.
  
\bibitem{Daans21}
Nicolas Daans, \emph{Universally defining finitely generated subrings of global
  fields}, Doc. Math. \textbf{26} (2021), pp.~1851--1869.

\bibitem{Daans22}
\bysame, \emph{Existential first-order definitions
and quadratic forms}, PhD Thesis in Mathematics defended at Antwerp University (2022).

\bibitem{Daans23}
\bysame, \emph{Universally defining $\mathbb{Z}$ in $\mathbb{Q}$ with 10
  quantifiers}, arXiv e-prints 2301.02107 (2023).

\bibitem{DDF}
Nicolas Daans, Philip Dittmann, and Arno Fehm, \emph{Existential rank and
  essential dimension of diophantine sets}, arXiv e-prints 2102.06941 (2021).

\bibitem{DPR}
Martin Davis, Hilary Putnam, and Julia Robinson, \emph{The decision problem for
  exponential diophantine equations}, Ann. of Math. (2) \textbf{74} (1961),
  pp.~425--436.

\bibitem{Dit18}
Philip Dittmann, \emph{Irreducibility of polynomials over global fields is
  diophantine}, Compos. Math. \textbf{154} (2018), no.~4, pp.~761--772.

\bibitem{DP}
Roberto Dvornicich and Laura Paladino, \emph{Local-global questions for
  divisibility in commutative algebraic groups}, Eur. J. Math. \textbf{8}
  (2022), pp.~S599--S628.

\bibitem{DZ}
Roberto Dvornicich and Umberto Zannier, \emph{Local-global divisibility of
  rational points in some commutative algebraic groups}, Bull. Soc. Math.
  France \textbf{129} (2001), no.~3, pp.~317--338.

\bibitem{DZ_04}
\bysame, \emph{An analogue for elliptic curves of the {G}runwald-{W}ang
  example}, C. R. Math. Acad. Sci. Paris \textbf{338} (2004), no.~1,
  pp.~47--50.

\bibitem{DZ2}
\bysame, \emph{On a local-global principle for the divisibility of a rational
  point by a positive integer}, Bull. Lond. Math. Soc. \textbf{39} (2007),
  no.~1, pp.~27--34.

\bibitem{Eis05}
Kirsten Eisentr\"{a}ger, \emph{Integrality at a prime for global fields and the
  perfect closure of global fields of characteristic {$p>2$}}, J. Number Theory
  \textbf{114} (2005), no.~1, pp.~170--181.

\bibitem{Eis07}
\bysame, \emph{Hilbert's tenth problem for function fields of varieties over number fields and {$p$}-adic fields}, J. Algebra \textbf{310} (2007), no.~2, pp.~775--792.

\bibitem{Eis08}
\bysame, \emph{Hilbert's tenth problem for function fields of characteristic zero}, Model theory with applications to algebra and analysis {V}ol.2, London Math. Soc. Lecture Note Ser., vol.~350, pp.~237--254, Cambridge University Press, Cambridge, 2008.

\bibitem{EM}
Kirsten Eisentr\"{a}ger and Travis Morrison, \emph{Universally and
  existentially definable subsets of global fields}, Math. Res. Lett.
  \textbf{25} (2018), no.~4, pp.~1173--1204.

\bibitem{ES}
Kirsten Eisentr\"{a}ger and Alexandra Shlapentokh, \emph{Hilbert's tenth problem over function fields of positive characteristic not containing the algebraic closure of a finite field}, J. Eur. Math. Soc. \textbf{19} (2017), no.~7, pp.~2103--2138.

\bibitem{Ershov}
Juri~L. Er\v{s}ov, \emph{On the elementary theory of maximal normed fields.},
  Dokl. Akad. Nauk SSSR \textbf{165} (1965), pp.~21--23.

\bibitem{Fehm}
Arno Fehm, \emph{Elementary geometric local-global principles for fields}, Ann.
  Pure Appl. Logic \textbf{164} (2013), no.~10, pp.~989--1008.

\bibitem{FJ}
Michael~D. Fried and Moshe Jarden, \emph{Field arithmetic}, third ed.,
  Ergebnisse der Mathematik und ihrer Grenzgebiete, 3. Folge, A Series of
  Modern Surveys in Mathematics, vol.~11, Springer-Verlag,
  Berlin, 2008, revised by Jarden.

\bibitem{GS06}
Philippe Gille and Tam\'{a}s Szamuely, \emph{Central simple algebras and
  {G}alois cohomology}, Cambridge Studies in Advanced Mathematics, vol. 165,
  Cambridge University Press, Cambridge, 2017.

\bibitem{GR}
Florence Gillibert and Gabriele Ranieri, \emph{On the local-global divisibility
  of torsion points on elliptic curves and {${\rm GL}_2$}-type varieties}, J.
  Number Theory \textbf{174} (2017), pp.~202--220.

\bibitem{Har_Sko}
Yonatan Harpaz and Alexei~N. Skorobogatov, \emph{Hasse principle for {K}ummer
  varieties}, Algebra Number Theory \textbf{10} (2016), no.~4, pp.~813--841.

\bibitem{Hodges}
Wilfrid Hodges, \emph{A shorter model theory}, Cambridge University Press,
  Cambridge, 1997.

\bibitem{Hsi_sil}
Liang-Chung Hsia and Joseph~H. Silverman, \emph{On a dynamical {B}rauer-{M}anin
  obstruction}, J. Th\'{e}or. Nombres Bordeaux \textbf{21} (2009), no.~1,
  pp.~235--250.

\bibitem{Kato}
Kazuya Kato, \emph{A {H}asse principle for two-dimensional global fields}, J.
  Reine Angew. Math. \textbf{366} (1986), pp.~142--183, with an appendix by
  Jean-Louis Colliot-Th\'{e}l\`{e}ne.

\bibitem{Katz}
Nicholas~M. Katz, \emph{Galois properties of torsion points on abelian
  varieties}, Invent. Math. \textbf{62} (1981), no.~3, pp.~481--502.


\bibitem{KimRoushC}
Ki-Hang Kim and Fred W.~Roush, \emph{Diophantine undecidability of {${\bf C}(t_1,t_2)$}}, J. Algebra \textbf{150} (1992), no.~1, pp.~35--44.

\bibitem{KimRoushpadic}
\bysame, \emph{Diophantine unsolvability over {$p$}-adic function fields}, J. Algebra \textbf{176} (1995), no.~1, pp.~83--110.

\bibitem{KoeSurvey}
Jochen Koenigsmann, \emph{Undecidability in number theory}, Model theory in algebra, analysis and arithmetic ({H}eidelberg, 2014), Lecture Notes in Math., 2111, 2014, pp.~159--195.
  
\bibitem{Koe16}
\bysame, \emph{Defining {$\Bbb Z$} in {$\Bbb Q$}}, Ann. of Math. (2)
  \textbf{183} (2016), no.~1, pp.~73--93.

\bibitem{Lam05}
Tsit~Y. Lam, \emph{Introduction to quadratic forms over fields}, Graduate
  Studies in Mathematics, vol.~67, American Mathematical Society, Providence,
  RI, 2005.

\bibitem{Malle}
Gunter Malle, \emph{Local-global conjectures in the representation theory of
  finite groups}, Representation theory---current trends and perspectives, EMS
  Ser. Congr. Rep., Eur. Math. Soc., Z\"{u}rich, 2017, pp.~519--539.

\bibitem{Mar02}
David Marker, \emph{Model theory}, Graduate Texts in Mathematics, vol. 217,
  Springer-Verlag, New York, 2002.

\bibitem{Mat}
Yuri~V. Matijasevi\v{c}, \emph{The {D}iophantineness of enumerable sets}, Dokl.
  Akad. Nauk SSSR \textbf{191} (1970), pp.~279--282.

\bibitem{MRFB}
Maruti {R}am {M}urty and {B}randon {F}odden, \emph{Hilbert's tenth problem}, An introduction to logic, number theory, and computability, Student Mathematical Library, vol.~88, American Mathematical Society, Providence, RI, 2019.

\bibitem{NSW}
J\"{u}rgen Neukirch, Alexander Schmidt, and Kay Wingberg, \emph{Cohomology of
  number fields}, second ed., Grundlehren der mathematischen Wissenschaften
  [Fundamental Principles of Mathematical Sciences], vol. 323, Springer-Verlag,
  Berlin, 2008.

\bibitem{Pal3}
Laura Paladino, \emph{On counterexamples to local-global divisibility in
  commutative algebraic groups}, Acta Arith. \textbf{148} (2011), no.~1,
  pp.~21--29.

\bibitem{PRV_12}
Laura Paladino, Gabriele Ranieri, and Evelina Viada, \emph{On local-global
  divisibility by {$p^n$} in elliptic curves}, Bull. Lond. Math. Soc.
  \textbf{44} (2012), no.~4, pp.~789--802.

\bibitem{PRV_14}
\bysame, \emph{On the minimal set for counterexamples to the local-global
  principle}, J. Algebra \textbf{415} (2014), pp.~290--304.

\bibitem{Par_Sur}
Raman Parimala and Venapally Suresh, \emph{Local-global principle for classical
  groups over function fields of {$p$}-adic curves}, Comment. Math. Helv.
  \textbf{97} (2022), no.~2, pp.~255--304.

\bibitem{Park}
Jennifer Park, \emph{A universal first-order formula defining the ring of
  integers in a number field}, Math. Res. Lett. \textbf{20} (2013), no.~5,
  pp.~961--980.

\bibitem{Phe91}
Thanases Pheidas,
\emph{Hilbert's tenth problem for fields of rational functions over finite fields}, Invent.~Math. \textbf{103} (1991), no.~1, pp.~1--8.

\bibitem{Poonen03}
Bjorn Poonen, \emph{Hilbert's tenth problem and {M}azur's conjecture for large subrings of {$\Bbb Q$}}, J. Amer. Math. Soc. \textbf{16} (2003), no.~4, pp.~981--990.

\bibitem{Poonen}
\bysame, \emph{Characterizing integers among rational numbers with a
  universal-existential formula}, Amer. J. Math. \textbf{131} (2009), no.~3,
  pp.~675--682.

\bibitem{Prestel}
Alexander Prestel, \emph{Pseudo real closed fields.}, Set theory and model
  theory ({B}onn, 1979), Lecture Notes in Math., 872, 1981, pp.~127--156.

\bibitem{robinson_1949}
Julia Robinson, \emph{Definability and decision problems in arithmetic}, J.
  Symbolic Logic \textbf{14} (1949), pp.~98--114.

\bibitem{robinson_1965}
\bysame, \emph{The decision problem for fields}, Theory of {M}odels ({P}roc.
  1963 {I}nternat. {S}ympos. {B}erkeley), North-Holland, Amsterdam, 1965,
  pp.~299--311.

\bibitem{Rum}
Robert~S. Rumely, \emph{Undecidability and definability for the theory of global fields}, Trans. Amer. Math. Soc. \textbf{262} (1980), no.~1, pp.~195--217.

\bibitem{Rum86}
Robert~S. Rumely, \emph{Arithmetic over the ring of all algebraic integers}, J. Reine Angew. Math. \textbf{368} (1986), pp.~127--133.

\bibitem{RS}
Russell Miller and Alexandra Shlapentokh, \emph{On existential definitions of c.e. subsets of rings of functions of characteristic 0}, Ann. Pure Appl. Logic \textbf{173} (2022), no.~4, Paper No.~103076, 50.

\bibitem{Sansuc}
Jean-Jacques Sansuc, \emph{Groupe de {B}rauer et arithm\'{e}tique des groupes
  alg\'{e}briques lin\'{e}aires sur un corps de nombres}, J. Reine Angew. Math.
  \textbf{327} (1981), pp.~12--80.

\bibitem{Ser_64}
Jean-Pierre Serre, \emph{Sur les groupes de congruence des vari\'{e}t\'{e}s
  ab\'{e}liennes.}, Izv. Akad. Nauk SSSR Ser. Mat. \textbf{28} (1964), pp.~3--20.

\bibitem{LocalFieldsSerre}
\bysame, \emph{Local fields}, Graduate Texts in Mathematics, vol.~67,
  Springer-Verlag, New York-Berlin, 1979, translated from the French by Marvin
  Jay Greenberg.

\bibitem{Shla_2000}
Alexandra Shlapentokh, \emph{Hilbert's tenth problem for algebraic function fields over infinite fields of constants of positive characteristic}, Pacific J. Math. \textbf{193} (2000), no.~2, pp.~463--500.

\bibitem{ShlBook}
\bysame, \emph{Hilbert's tenth problem}, Diophantine classes and extensions to global fields, New Mathematical Monographs,  vol.~7, Cambridge University Press, Cambridge, 2007.

\bibitem{Shla_2008}
\bysame, \emph{Elliptic curves retaining their rank in finite extensions and
  {H}ilbert's tenth problem for rings of algebraic numbers}, Trans. Amer. Math. Soc. \textbf{360} (2008), no.~7, pp.~3541--3555.

\bibitem{Shla_2012}
\bysame, \emph{Elliptic curve points and {D}iophantine models of {$\Bbb Z$} in large
  subrings of number fields}, Int. J. Number Theory \textbf{8} (2012), no.~6, pp.~1335--1365.

\bibitem{AEC}
Joseph~H. Silverman, \emph{The arithmetic of elliptic curves}, second ed.,
  Graduate Texts in Mathematics, vol. 106, Springer, Dordrecht, 2009.

\bibitem{Sko}
Alexei~N. Skorobogatov, \emph{Torsors and rational points}, Cambridge Tracts in
  Mathematics, vol. 144, Cambridge University Press, Cambridge, 2001.

\bibitem{Sus}
Andrei~A. Suslin, \emph{The structure of the special linear group over rings of
  polynomials.}, Izv. Akad. Nauk SSSR Ser. Mat. \textbf{11} (1977), no.~2, pp.~235--252,
  477.

\bibitem{Sut}
Andrew~V. Sutherland, \emph{A local-global principle for rational isogenies of
  prime degree}, J. Th\'{e}or. Nombres Bordeaux \textbf{24} (2012), no.~2,
  pp.~475--485.

\bibitem{Tarski}
Alfred Tarski, \emph{A {D}ecision {M}ethod for {E}lementary {A}lgebra and
  {G}eometry}, The Rand Corporation, Santa Monica, Calif., 1948.

\bibitem{Vogt}
Isabel Vogt, \emph{A local-global principle for isogenies of composite degree},
  Proc. Lond. Math. Soc. (3) \textbf{121} (2020), no.~6, pp.~1496--1530.

\bibitem{Quats}
John Voight, \emph{Quaternion algebras}, Graduate Texts in Mathematics, vol.
  288, Springer, Cham, 2021.

\bibitem{Wat}
William~C. Waterhouse, \emph{Abelian varieties over finite fields}, Ann. Sci.
  \'{E}cole Norm. Sup. (4) \textbf{2} (1969), pp.~521--560.

\bibitem{ZS}
Geng-Rui Zhang and Zhi-Wei Sun, \emph{{$\Bbb Q\setminus\Bbb Z$} is diophantine
  over {$\Bbb Q$} with 32 unknowns}, Bull. Pol. Acad. Sci. Math. \textbf{70}
  (2022), no.~2, pp.~93--106.

\end{thebibliography}

\providecommand{\bysame}{\leavevmode\hbox to3em{\hrulefill}\thinspace}
\providecommand{\MR}{\relax\ifhmode\unskip\space\fi MR }
\providecommand{\MRhref}[2]{%
  \href{http://www.ams.org/mathscinet-getitem?mr=#1}{#2}
}
\providecommand{\href}[2]{#2}

\end{document}